\documentclass[12pt]{article}
\usepackage{amsmath,amsthm,amsfonts,graphicx,amssymb}
\usepackage{url}
\newtheorem{Theorem}{Theorem}[section]
\newtheorem{lemma}[Theorem]{Lemma}
\newtheorem{corollary}[Theorem]{Corollary}

\theoremstyle{definition}
\newtheorem{definition}[Theorem]{Definition}

\newtheorem{remark}[Theorem]{Remark}

\newtheorem{example}[Theorem]{Example}
\numberwithin{equation}{section}

\def\ep{\hfill $\Box$}
\def\epcl{\hfill $\diamondsuit $}
\def\bp{\noindent{\bf Proof}  \ }

\begin{document}

\title{{\bf Spectra of Digraph Transformations}
\thanks{The research is supported in part by the Fundamental Research Funds for the Central Universities of China
11D10902. }}

\author{ {\bf Aiping Deng}$^{\rm a,}$\footnote{Corresponding author.
Email: apdeng@dhu.edu.cn. Tel: 86-21-67792089-568. }\ ,
{\bf  Alexander Kelmans}$^{\rm b,c}$ \\
$^{\rm a}${\small\em{Department of Applied Mathematics, Donghua University, 201620 Shanghai, China} } \\
$^{\rm b}${\small\em{Department of Mathematics, University of Puerto Rico, San Juan, PR, United States}}\\
$^{\rm c}${\small\em{RUTCOR, Rutgers University, New Brunswick, NJ, United States}}
  }

\date{}
\maketitle

\begin{abstract}

Let $D = (V,E)$ be a  directed graph (or a digraph) with the vertex set $V = V(D)$ and the arc set 
$E = E(D) \subseteq V \times V \setminus \{(v,v): v \in V\}$ (and so $D$ has no loops and no multiple arcs).
%
%
Let $D^0$ be the digraph with vertex set
$V$ and with no arcs, $D^1$  the
complete digraph
with vertex set $V$,
$D^+ = D$, and $D^-$ the complement
$D^c$ of $D$.
For $e = (x, y) \in E$ let $x = t(e)$ and $y = h(e)$.
Let ${\cal T}(D)$ (${\cal T}^{cb}(D)$) be the digraph with vertex set $V\cup E$
such that
$(v, e)$ is an arc in
${\cal T}(D)$ (resp., in ${\cal T}^{cb}(D)$) if and only if $v \in V$, $e \in E$, and  vertex $v = t(e)$  (resp., $v \ne t(e)$)  in $D$.
Similarly, let ${\cal H}(D)$ (${\cal H}^{cb}(D)$) be the digraph with vertex set $V\cup E$
such that
$(e,v)$ is an arc in
${\cal H}(D)$ (resp., in ${\cal H}^{cb}(D)$) if and only if $v \in V$, $e \in E$, and  vertex $v = h(e)$
(resp., $v \ne h(e)$)  in $D$.
Given a digraph $D$ and three variables
$x, y, z \in \{0,1, +, -\}$, the {\em ${xyz}$-transformation} of $D$ is the digraph $D^{xyz}$
such that
$D^{xy0} = D^x\cup (D^l)^y$
and $D^{xyz} = D^{xy0}\cup W$, where
$W = {\cal T}(D) \cup {\cal H}(D)$  if $z = +$,
$W = {\cal T}^{cb}(D) \cup {\cal H}^{cb}(D)$ if $z = -$, and
$W$ is the complete bipartite digraph with parts $V$ and $E$ if $z = 1$.
In this paper we obtain the adjacency characteristic polynomials of
some $xyz$-transformations of an $r$-regular digraph $D$  in terms of the adjacency polynomial, the  number of vertices  of $D$ and $r$.
Similar  results are obtained for some non-regular digraphs, named {\em digraph-functions}.
Using $xyz$-transformations we   give various constructions of non-isomorphic cospectral digraphs.
Our notion of $xyz$-transformation and the corresponding adjacency polynomials results are also valid for digraphs with loops and multiple arcs provided
$x, y, z \in \{0, +\}$ and $z \in \{0, 1, +, -\}$.
We also extend the notion of
$xyz$-transformation and the above adjacency polynomial results to digraphs $(V,E)$ with possible 
loops and no multiple arcs.
\\[1ex]
\\[1ex]
\end{abstract}

\section{Introduction}

\indent

We will  consider finite digraphs with possible
loops and multiple arcs.
%
%
All notions on graphs and matrices that are used but not defined here can be found in
\cite{B&M07,D,G,horn,west}.

Let ${\cal D}$ denote the set of digraphs.
Various important results in graph theory have been obtained by considering some functions
$F: {\cal D} \to  {\cal D}$
or  $F_s: {\cal D}_1\times  \ldots \times {\cal D}_s \to  {\cal D}$ called
{\em operations}
(here  each ${\cal D}_i = {\cal D}$) and
by establishing how these operations affect certain properties or parameters of graphs or digraphs.
The complement,  the $k$-th power of a (di)graph, and  the line (di)graph are well known  examples of such operations.
Also, the  Bondy-Chvatal and Ryz\'{a}\u{c}ek closers of graphs are very useful operations in graph Hamiltonicity theory
\cite{B&M07}.
(Strengthenings and extensions of the
Ryz\'{a}\u{c}ek result are given in \cite{Kclosure}). Graph operations introduced by Kelmans in \cite{Ktrlp,KoperR}
turned out to be very  useful  because  they are  monotone with respect to some partial order relations on the set of graphs \cite{Kmxtr,Kprobcmpr}.
Gross and Tucker introduced the operation of voltage lifting on a graph which can be generalized to digraphs \cite{D&W05,G&T87}. By this operation one can obtain the derived covering (di)graph and the relation between the adjacency characteristic polynomials of the
(di)graph and its derived covering (di)graph \cite{DS&W07,D&W05,M&S95}.

The goal of  this paper is to  consider (and establish some properties of) certain operations depending on parameters  $x, y, z \in \{0,1,+,-\}$.
These operations  induce  functions $T^{xyz}: {\cal D} \to {\cal D}$. We put  $T^{xyz}(D) = D^{xyz}$ and call
$D^{xyz}$ the {\em $xyz$-transformation of $D$}, which is similar to the ${xyz}$-transformation of an undirected graph
(see, for example, \cite{DKM}).

For an undirected graph $G$, some graph properties
of the transformations $G^{xyz}$ with
$x,y,z\in \{+,-\}$ were discussed in
\cite{linshu, wubaoy2, wumeng1}.
For a regular undirected graph $G$, the adjacency polynomials  and  spectra of $G^{00+}$,
$G^{+0+}$, $G^{0++}$, and $G^{+++}$ were given in \cite{CDS} (pages 63 and 64).
Yan and Xu obtained the adjacency spectra of the other seven transformations $G^{xyz}$ with $x, y, z \in \{ +, - \}$ in terms of the adjacency spectrum of $G$
\cite{yanxu}.
In 1967 Kelmans established the formulas for the Laplacian polynomials and the number of spanning trees of $G^{0++}, G^{0+0}, G^{00+}$, and $L(G)$ \cite{kelmans}. Recently, Deng, Kelmans, and Meng presented for a regular graph $G$ and all $x, y, z \in \{0,1,+, -\}$ the formulas
of the Laplacian  polynomials and the number of spanning trees of $G^{xyz}$ in terms of
the number of vertices, number of edges, and the Laplacian spectrum of $G$ \cite{DKM}.
The zeta functions of $G^{0++}$ and $G^{+++}$ and their coverings were discussed in \cite{KwakSato}.
The transformations $G^{00+}$,
$G^{0++}$, and $G^{+++}$ for a (di)graph $G$ were also called {\em subdivision $($di$)$graph, middle $($di$)$graph}, and {\em total $($di$)$graph of $G$}, respectively  \cite{CDS,kelmans,KwakSato, zhang}.

In \cite{zhang} Zhang, Lin and Meng presented the  adjacency polynomials of $D^{00+}$, $D^{+0+}$, $D^{0++}$, and  $D^{+++}$ for  any digraph $D$. The adjacency polynomials (and spectra) of the other transformations $D^{xyz}$ of
a regular digraph $D$ with $x, y, z \in \{ +, - \}$ were obtained by Liu and Meng \cite{liumeng}.

In this paper we
give  descriptions of
the  adjacency characteristic polynomials of $xyz$-transformation of all $r$-regular digraphs for all  $x, y, z \in \{0,1,+,-\}$
  as well as for some non-regular digraphs, for example, for so-called digraph-functions.
From these descriptions it follows that
the spectrum $S_a(D^{xyz})$ of the $xyz$-transformation of any $r$-regular digraph $D$ is uniquely defined by
the spectrum $S_a(D)$ of $D$.
Moreover, we  obtain the explicit description of $S_a(D^{xyz})$ in terms of $S_a(D)$, $r$, and $v(D)$, the number of vertices of $D$ (where $r$ and $v(D)$ are uniquely defined by $S_a(D)$).
The results of this paper may be considered as a natural and useful extension of the results in  \cite{CDS},
Section 2  ``Operations on Graphs and the Resulting  Spectra''.

In Section \ref{notions} we introduce main
notions, notation, and simple observations.
 Some preliminaries are given in Section \ref{preliminaries}.
 In Section \ref{z-in-{0,1}} we
  present the adjacency polynomials  of some transformations $D^{xyz}$ with $z \in \{0,1\}$.
In Section \ref{z-in-{+,-}} we describe the adjacency polynomials of some transformations
$D^{xyz}$ with $z\in \{+,-\}$ and $\{x,y\} \cap \{0,1\} \neq \emptyset$ for regular digraphs.
We also were able to obtain similar results for some non-regular digraphs.
In Section \ref{digraph-functions} we
consider a special class of non-regular digraphs called the digraph-functions, give some criteria for a digraph to be a digraph-function,  and
describe the adjacency polynomials of some transformations
$D^{xyz}$ for  digraph-functions $D$ and their inverse.
In Section \ref{cospectral} we summarize some previous constructions providing various pairs of non-isomorphic and cospectral digraphs (including $xyz$-transformations) and give some more results of this nature.
Section \ref{remarks} contains some additional remarks and questions.
In Appendix  we provide for all $x, y, z \in \{0,1,+, -\}$ the list of formulas
for the adjacency polynomials of the $xyz$-transformations of an $r$-regular digraph $D$ in terms of $r$,
the number of vertices $n$, the number of edges $m = nr$, and
the adjacency polynomial of $D$.

\section{Some notions and notation}

\label{notions}

\indent

A {\em general directed graph} $D$
(or {\em a digraph with possible multiple arcs and loops})
is a triple  $(V, E, \psi )$, where $V$ and
$E$ are finite sets, $V$ is non-empty,
and $\psi $ is a function from $E$ to
$V\times V$ (and so
 $\psi (e)$ is the ordered pair of {\em ends} of arc  $e$ in  $E$). If $\psi (e) = (v,v)$ for some
$v \in V$, then arc $e$ is called a {\em loop} in $D$.
The sets $V$ and $E$ are called the {\em vertex set } and the {\em arc set } of digraph $D$ and denoted by $V(D)$ and $E(D)$, respectively.
Let $v(D) = |V(D)|$ and $e(D) = |E(D)|$.

Given two general digraphs
$D_1 = (V_1, E_1, \psi _1)$ and
$D_2 = (V_2, E_2, \psi _2)$, a pair
$(\alpha _v, \alpha _e)$ of bijections
$\alpha _v: V_1 \to V_2$ and
$\alpha _e: E_1 \to E_2$ is called an {\em isomorphism from $D_1$ to $D_2$} if for every $a \in E_1$,
$ \psi _1(a) = (x,y) \Leftrightarrow
 \psi _2(\alpha _e (a)) =
 (\alpha _v(x), \alpha _v(y))$.
 We say that {\em digraph $D$ is isomorphic to digraph  $F$} (or equivalently, {\em $D$ and $F$ are isomorphic}) and write  $D \cong F$ if there exists an isomorphism from $D$ to $F$.

A {\em  directed graph} (or a {\em digraph}) is a general digraph $(V, E, \psi )$, where function
$\psi : E \to V\times V$ is injective.
In other words, a digraph
$D$ is a pair $(V, E)$, where $V$ is a non-empty set and $E \subseteq V\times V$, and so $D$ has no multiple arcs but may have at most one loop in each vertex.
If $E = V\times V$, then digraph $D$ is called
a {\em complete digraph} and denoted by $K_{\circ}$.

Given two digraphs
$D_1 = (V_1, E_1)$ and
$D_2 = (V_2, E_2)$, a  bijection
$\alpha : V_1 \to V_2$  is called an {\em isomorphism from $D_1$ to $D_2$} if
$(x,y) \in E_1 \Leftrightarrow
  (\alpha (x), \alpha (y)) \in E_2$.
 As above, we say that {\em digraph $D$ is isomorphic to digraph  $F$} (or equivalently, {\em $D$ and $F$ are isomorphic}) and write  $D \cong F$ if there exists an isomorphism from $D$ to $F$.

A digraph $D = (V, E)$ is called {\em simple} if $D$ has  no loops, and so $E \subseteq
\{V\times V\}$, where
$\{V\times V\} = V\times V \setminus \{(x,x): x \in V\}$.
Given a  digraph $D$ with $V = V(D) = V(K_\circ )$, let
$D^c_\circ = K_\circ \setminus E(D)$. Digraph
 $D^c_\circ $ is called the {\em $K_\circ $-complement of} $D$.
Let $K$ be the graph obtained from $K_\circ$ by removing all its loops, i.e. $E(K) = \{V\times V\}$. We call $K$
a {\em simple complete digraph}.
Given a simple digraph $D$, let
$D^c = K \setminus E(D)$.
Digraph $D^c$ is called the
{\em $K$-complement } (or simply, {\em complement}) {\em of }$D$.

A digraph $D = (V, E)$ is called {\em $(X,Y)$-bipartite} if
$V = X\cup Y$, $X \cap Y = \emptyset$, and $E \subseteq  X \times Y$. If, in addition, $E = X \times Y$,
then $D$ is called a {\em complete $(X,Y)$-bipartite digraph} and is denoted by $K_{XY}$.
Given an $(X,Y)$-bipartite digraph, let
$D^{cb} = K_{XY} \setminus E(D)$. Digraph $D^{cb}$ is called the {\em $(X,Y)$-bipartite complement of} $D$.

If $e \in E$ is an arc in $D$ and
$\psi (e) = (u,v)$
(possibly, $u = v$), then $u$ is called the {\em tail} of arc $e$ and $v$  called the {\em head} of arc $e$ and we put
$t(e,D) = t(e) =  u$ and $h(e,D) = h(e) = v$.
Obviously, $t$ and $h$ are functions:
$t: E \to V$ and  $h: E \to V$.
The functions $t$ and $h$ can also be described by the corresponding $(V \times E)$-matrices $T$ and $H$:
\\
$(t)$
the \textit{tail incidence matrix} $T(D) = T = \{t_{ij}\}$ of $D$, where $t_{ij} = 1$ if $v_i = t(e_j)$ and
$t_{ij} = 0$, otherwise,
and
\\
$(h)$
the \textit{head incidence matrix}
$H(D) = H = \{h_{ij}\}$ of $D$,
where $h_{ij}= 1$ if $v_i= h(e_j)$ and $h_{ij}  = 0$, otherwise.

The {\em line digraph} of $D$, denoted by $D^l$, is a digraph with vertex set $E(D)$ and arc set
$E(D^l)=\{(p, q) : p, q \in E(D) \mbox{\ and \ }
h(p,D) = t(q, D)\}$.
Obviously, $D^l$ has no multiple arcs and if $D$ is simple, then $D^l$ is also simple.

For $v \in V(D)$, let
$d_{out}(v, D) = d_{out}(v) = |\{e \in E(D): t(e) = v\}|$
be the {\em out-degree of} $v$ and
$d_{in}(v, D) = d_{in}(v) = |\{e \in E(D): h(e) = v\}|$ be
the {\em in-degree of $v$} in $D$.

A digraph $D$ is called {\em balanced} if
$d_{in}(v, D) = d_{out}(v, D) \ne 0$ for every $v \in V(D)$, and
$D$ is called $r$-{\em regular} if
$d_{in}(v, D) = d_{out}(v, D)  = r $ for every $v \in V(D)$, and so every $r$-regular digraph with $r \ge 1$ is balanced.

An undirected graph $G$ is the {\em underlying graph} of a digraph $D$ if
$V(G) = V(D)$ and $[x,y] \in E(G)$ if and only if $x \ne y$ and either $(x,y) \in E(D)$ or $(y,x) \in E(D)$.
 A digraph $D$ is {\em connected} if  its underlying graph is connected and
{\em not connected}, otherwise. A {\em component} of a digraph $D$ is a maximal connected subdigraph of $D$.
Obviously, two different components of $D$ are disjoint
(i.e. have no common vertex).

A digraph $D$ is called {\em strongly connected} if $D$ has a directed path from  $x$ to  $y$ for every
ordered pair $(x,y)$ of vertices in $D$.
Obviously, a connected balanced digraph is strongly connected.

For a digraph $D = (V, E)$, let
$D^{-1} = (V, E^{-1})$, where  $E^{-1} = \{(y,x): (x,y) \in E\}$. The digraph $D^{-1}$ is called the {\em inverse of} $D$.
Let $D^0$ be the digraph with vertex set $V(D)$ and with no arcs, $D^1$  the simple complete digraph with vertex set $V(D)$,
$D^+ = D$ and $D^- = D^c$ if $D$ is simple.

Let ${\cal T}(D)$ (${\cal T}^{cb}(D)$) denote the digraph with vertex set $V \cup E$ and such that
$(v, e)$ is an arc in
${\cal T}(D)$ (resp., in ${\cal T}^{cb}(D)$) if and only if $v \in V$, $e \in E$, and  vertex $v = t(e)$  (resp., $v \ne t(e)$)  in $D$, and so ${\cal T}(D)$ is a $(V,E)$-bipartite digraph and ${\cal T}^{cb}(D)$ is the $(V,E)$-bipartite complement of ${\cal T}(D)$.
Similarly, let ${\cal H}(D)$ (${\cal H}^{cb}(D)$) be the digraph with vertex set $V\cup E$ and such that
$(e,v)$ is an arc in
${\cal H}(D)$ (resp., in ${\cal H}^{cb}(D)$) if and only if $v \in V$, $e \in E$, and  vertex $v = h(e)$
(resp., $v \ne h(e)$)  in $D$, and so ${\cal H}(D)$ is an $(E,V)$-bipartite digraph and ${\cal H}^{cb}(D)$ is the
$(E,V)$-bipartite complement of ${\cal H}(D)$.

Given two digraphs $D$ and $D'$, let
$D \cup D'$ denote the digraph with
$V(D \cup D') = V(D) \cup V(D')$ and
$E(D \cup D') = E(D) \cup E(D')$.

The main notion of the digraph transformations we are  going to discuss is using the notion of
$K$-complement of a digraph $D$ and is defined as follows.
\begin{definition}
\label{definition}
Given a simple digraph $D$ and three variables
$x, y, z \in \{0,1, +, -\}$, the {\em ${xyz}$-transformation $D^{xyz}$  of} $D$ is the digraph
such that
$D^{xy0} = D^x\cup (D^l)^y$
and $D^{xyz} = D^{xy0}\cup W$, where
$W = {\cal T}(D) \cup {\cal H}(D)$  if $z = +$,
$W = {\cal T}^{cb}(D) \cup {\cal H}^{cb}(D)$ if $z = -$, and
$W$ is the union of complete
$(V, E)$-bipartite and $(E,V)$-bipartite digraphs if $z = 1$.
\end{definition}

From the definitions of ${\cal T}(D)$ and ${\cal H}(D)$ we have:
\begin{remark}
\label{subdividing}
{\em Digraph
${\cal T}(D) \cup {\cal H}(D)$ can be obtained from $D$ by subdividing every arc $e$ of $D$
into two arcs by a new vertex with label $e$, and it is just the subdivision digraph $D^{00+}$ of $D$.}
\end{remark}

\begin{remark}
\label{general-digraph-transform}
{\em Definition \ref{definition} is  also valid
for general digraphs $D$ in the case when  $x , y \in \{0,+\}$
and $z \in \{0,1,+,-\}$.}
\end{remark}

Obviously, if $D$ is a simple digraph, then  $D^{xyz}$ is also a simple digraph for $x,y,z \in \{0,1,+,-\}$.

\begin{figure}[ht]
\begin{center}
\scalebox{0.7}[.7]{\includegraphics{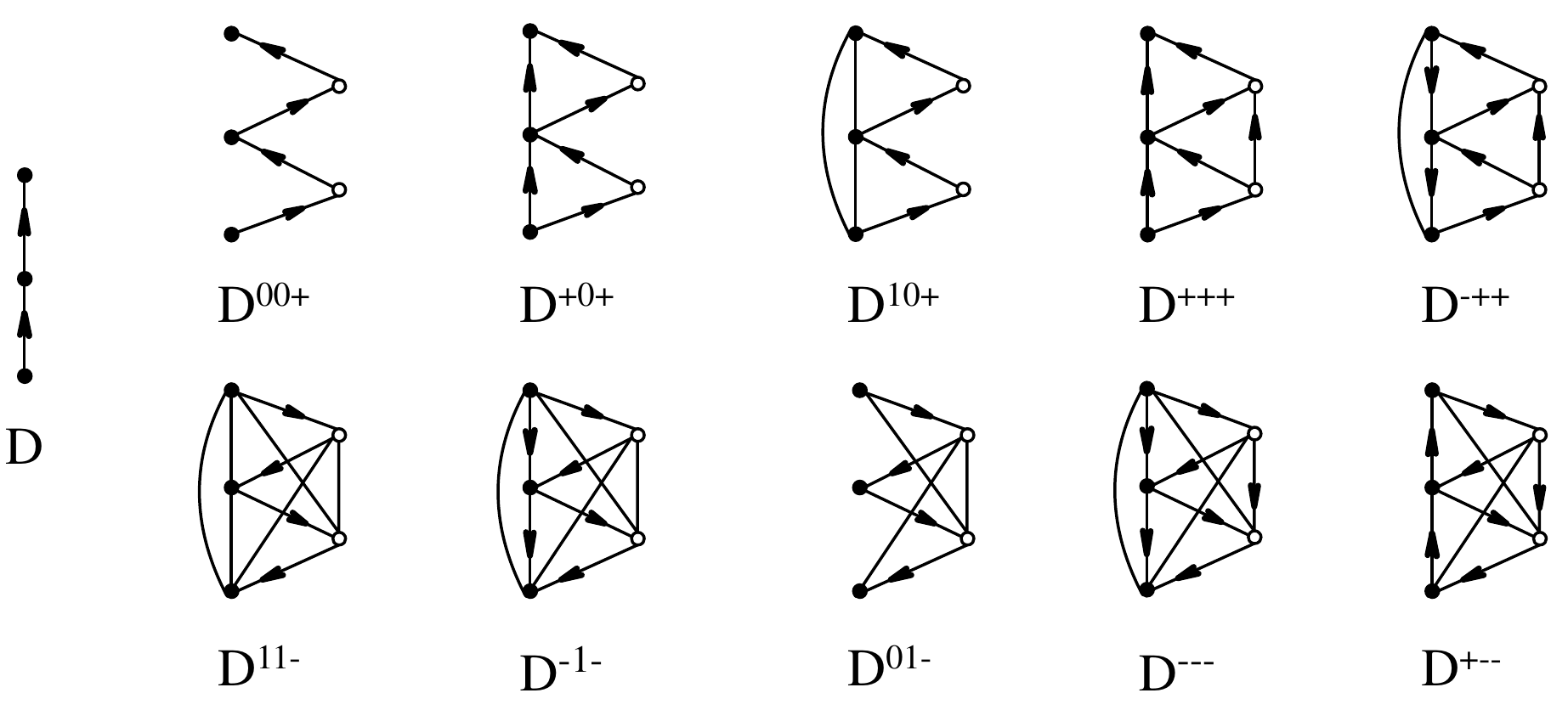}}
\end{center}
\caption{Digraph $D$ and some of its $xyz$-transformations. The undirected edge represents a pair of arcs having the opposite directions.}
\label{trDi}
\end{figure}

Examples of the  ${xyz}$-transformations of a 3-vertex directed path $D$ are given in Figure \ref{trDi}.

Let $V = V(D) = \{v_1, \ldots , v_n\}$ and
$E = E(D)  = \{e_1, \ldots , e_m\}$.
Let $A(D)$ be the $(V\times V)$-matrix $(a_{ij})$ such that
$a_{ij} = |\{e \in E(D) : \psi (e) = (v_i, v_j)\}|$
and $R(D)$  the (diagonal) $(V\times V)$-matrix $(r_{ij})$ such that $r_{ii} = d_{out}(v_i)$ and $r_{ij} = 0$ for $i \ne j$.

Now let $L(D) = R(D) - A(D)$.
The matrices $A(D)$ and $L(D)$ are called the  {\em adjacency} and the {\em Laplacian matrix} of $D$, respectively.

Let $I_n$ be the identity $(n\times n)$-matrix and
$J_{m n}$ the all-ones $(m\times n)$-matrix.
Obviously, if $v(D) = n$, then
$A(D^c) = J_{nn} - I_n - A(D)$ and
$A(D^c_{\circ }) = J_{nn}  - A(D)$.

The characteristic polynomials
$A(\lambda, D) = \det (\lambda I - A(D))$ and
$L(\lambda, D) = \det (\lambda I - L(D))$ of matrices $A(D)$
and $L(D)$ are called {\em the adjacency}  and
the {\em Laplacian polynomials of} $D$, respectively, and the sets $S_a(D)$ and $S_l(D)$ of roots of
$A(\lambda, D)$ and of $L(\lambda, D)$
(with their multiplicities) are the  corresponding spectra of $D$.

We call digraphs $D$ and $F$
{\em adjacency cospectral} or simply {\em cospectral} and write $D \sim ^A F$ if
$A(\lambda, D) = A(\lambda, F)$.

In what follows, we will often omit $D$ from the notation
by writing $V$ instead of $V(D)$, $E$ instead of $E(D)$, $A$ instead of $A(D)$, etc.  when the digraph $D$ is clear from the context.

\section{Preliminaries}
\label{preliminaries}

\indent

We start with some simple and useful observations about digraph $D^{xyz}$.
\begin{lemma}
\label{Gxyz}
Let $D= (V,E)$ be a simple digraph and $x,y,z \in
\{0, 1, +, -\}$. Then
\\[1ex]
$(a1)$ graphs $D^{xyz}$ and $D^{x'y'z'}$ are complement if and only if each of $\{x,x'\}$, $\{y,y'\}$, $\{z,z'\}$  is either
$\{0,1\}$ or  $\{+,-\}$
and
\\[1ex]
$(a2)$ if $K$ is a simple complete digraph, then
$K^{0yz} =  K^{-yz}$ and  $K^{x0z} = K^{x-z}$ as well as
$K^{1yz} =  K^{+yz}$ and
$K^{x1z} = K^{x+z}$.
\end{lemma}

\begin{lemma}
\label{Inverse}
Let $D$ be a digraph and $D^{-1}$ be the inverse of $D$. Then
\\[1ex]
$(a1)$
$A(\lambda , D) = A(\lambda , D^{-1})$ and
\\[1ex]
$(a2)$
$A(\lambda , (D^{-1})^{xyz}) =
A(\lambda , (D^{xyz})^{-1}) = A(\lambda , (D^{xyz}))$ for
$x,y,z \in \{0, 1, +, -\}$.

\end{lemma}

Using Theorem 1.2 in \cite{CDS} it is easy to prove the following spectrum property  of line digraphs.
\begin{lemma}
\label{PL}
Let $D$ be a digraph with $n$ vertices and $m$ arcs. Then
$$A(\lambda, D^l) = \lambda^{m-n} A(\lambda, D).$$
\end{lemma}

We will need the following two simple lemmas
on the matrices $A = A(D)$, $A^l = A(D^l)$, $H = H(D)$, and $T = T(D)$.
\begin{lemma}
\label{AL}
Let $D$ be a digraph. Then
\\[0.2ex]
$(a1)$ $A = T H^\top $ and
\\[0.2ex]
$(a2)$ $A^l = H^\top T $.
\end{lemma}

\begin{lemma}
\label{lemJ}
Let $D$ be a simple  $r$-regular  digraph with $n$ vertices and $m$ arcs.
Let $k$ be a positive integer. Then
\\[1ex]
$(a1)$
$T J_{m k} = r J_{n k}$,
\\[1ex]
$(a2)$
$J_{k n} T =  J_{k m}$,
\\[1ex]
$(a3)$
$J_{k m} H^{\top} = r J_{k n}$,
\\[1ex]
$(a4)$
$H^{\top} J_{n k}=  J_{m k}$,
\\[1ex]
$(a5)$
$J_{k n}A = r J_{k n}$, and
\\[1ex]
$(a6)$
$A J_{n k} = r J_{n k}$.
\end{lemma}

We will also need the following classical  fact on matrices.
\begin{lemma}
\label{lemABCD}
{\em \cite{G,horn}}
Let $A$ and $D$ be square matrices. Then
 \[ \left|\begin{array}{cc}
A &\quad B\\ C &\quad D \end{array}\right|\quad
=\left\{\begin{array}{ll}|A|~|D - C A^{-1} B|, &\quad \,\, if\, $A$
\mbox{\, is\, invertible},
\\[1ex]
|D|~|A - B D^{-1}C|, &\quad \,\, if\, $D$
\mbox{\, is\, invertible}.
\end{array}\right.\]
\end{lemma}

The other preliminaries we give below include the important Reciprocity Theorem on the relation between the Laplacian spectra of the complement digraphs and the corresponding Reciprocity Theorem for the adjacency spectra of regular digraphs $D$ and $D^c$
\cite{Ktree1}. In particular, because of this theorem  it is sufficient to describe the   adjacency characteristic polynomials of $xyz$-transformations of regular digraphs up to
the graph operation of taking the complement.
\begin{lemma} {\em \cite{Kcourse,KelLatvia}}
\label{multipl=cmp}
Let $D$ be a digraph such that each component of $D$
is strongly connected. Then the multiplicity of the zero eigenvalue
of $L(D)$ is equal to the number of components of $D$.
\end{lemma}

\begin{lemma} {\em \cite{Kcourse,KelLatvia}}
\label{ortogLaplac}
Let $D$ be an balanced digraph with $n$ vertices.
Then $L = L(D)$ has a set of   eigenvectors
$X_1,  \cdots , X_n$ such that $X_n = J_{n1}$ and $X_i$ is  orthogonal to $X_n$ for every
$i =1, 2, \ldots, n-1$.
\end{lemma}
\bp
Obviously, every connected  balanced digraph is strongly connected, and so every component of $D$ is strongly connected.

By definition of $L = L(D)$, $LJ_{n 1} = 0$, and so $J_{n1}$ is  an eigenvector of $L$ corresponding to a zero eigenvalue. Since $D$ is balanced, we also have: $J_{1 n} L= 0$.
Let $V(D) = V$. We can interpret every eigenvector $X_i$ of $D$ as a function from $V$ to
$ \mathbb{R}$.
\\[1ex]
{\bf (p1)} Suppose that $D$ is connected. Then by Lemma
\ref{multipl=cmp}, the multiplicity of the zero eigenvalue of $L$ is equal to one.
Let $\lambda _i$, $i = 1 \ldots , n$, be an eigenvalue and  $X_i$ be the corresponding eigenvector of $L$ and let $\lambda _n = 0$ and $X_n = J_{n1}$. Then $\lambda _i \ne 0$ for $i \ne n$.
Now
$0 X_i = J_{1n}L X_i = J_{1n} \lambda _i X_i =
\lambda _i J_{1n} X_i$, $i \ne n$.
Therefore $J_{1n} X_i = 0$ for $i \ne n$, i.e. each
$X_i$ is orthogonal to $X_n = J_{n1}$.
\\[1ex]
{\bf (p2)} Now suppose that $D$ is not connected.
Let $D_1, \cdots , D_k$ be components of $D$, and so each $D_s$ is strongly connected. Then by
Lemma \ref{multipl=cmp}, the multiplicity of the zero eigenvalue of $L$ is equal to $k$.
Obviously,
$L(\lambda , D) = \prod _{s = 1}^k L(\lambda , D_s)$, and so
$S_l(D) = \cup_{s = 1}^k S_l(D_s)$.
Let $n_s$ be the number of vertices of $D_s$, and so
$n_1 + \cdots + n_k = n$.
Let $\lambda ^s_i$ be an eigenvalue and  $x^s_i$ the corresponding eigenvector of $L(D_s)$.
Set $\lambda^s _{n_s} = 0$ and then $x^s_{n_s} = J_{n_s1}$. Let $V_s = V(D_s)$. Then
$x^s_i$ can be interpreted as a function from $V_s$ to
$ \mathbb{R}^{n_s}$. Let $X^s_i$ be a function from $V$ to $ \mathbb{R}^n$ such that $X^s_i(v) = x^s_i(v)$ for $v \in V_s$ and
$X^s_i(v) = 0$ for $v \in V \setminus V_s$.
By ${\bf (p1)}$, $\lambda^s _{n_s}$ is the only zero eigenvalue of $L(D_s)$ and
$J_{1n_s}x^s_i = 0$
for $i \ne n_s$. Therefore
$L X^s_i = \lambda ^s_iX^s_i $ and $J_{1n}X^s_i = 0$, and so
$X^s_i$ is an eigenvector of $L$ corresponding to its eigenvalue $\lambda ^s_i$ and $X^s_i$ is orthogonal to $X_n = J_{n1}$.

Now consider $Y_s = X^s_{n_s}$.
Then each $L Y_s = 0Y_s$, i.e.,
$Y_s$ is an eigenvector of $L$ corresponding to a zero eigenvalue of $L$. Moreover, $\{Y_1, \cdots , Y_k\}$ is a basis of the $k$-dimensional eigenspace $Q \subset  \mathbb{R}^n$ of $L$ corresponding to the zero eigenvalues, and
$X_n = J_{n1} = Y_1 + \cdots + Y_k \in Q$.
Then each $Z \in Q$ is a linear combination of $Y_j$'s: $Z = \delta _1 Y_1+ \cdots + \delta  _kY_k$, where each $\delta _s \in \mathbb{R}$. Obviously, $J_{1n} Z = 0$ if and only if $ \delta  _1n_1 + \cdots + \delta  _k n_k = 0$.
This condition defines the $k-1$ dimensional subspace
$Q'$ in $Q$ orthogonal to $X_n$. Let
$\{Z_1, \cdots , Z_{k-1}\}$ be a basis of $Q'$. Then
each $Z_j$ is an eigenvector of $L$ corresponding to
a zero eigenvalue and $Z_j$ is orthogonal to $X_n$.
\ep
\\[2ex]
\indent
Here is the important Reciprocity Theorem for  the Laplacian spectrum of simple balanced digraphs.
Let $S'_l(D)$ be the set of all Laplacian eigenvalues of $D$ except for one zero eigenvalue.
\begin{Theorem}
{\em \cite{KelLatvia,Ktree1}}
\label{di-reciprocity}
Let  $D$ be a  simple balanced digraph with $n$ vertices.
Then
\\[1ex]
$(a1)$ there exists a bijection
$\sigma :
S'_l(D) \to S'_l(D^c)$ such that $x + \sigma ( x) = n$
for every $x \in S'_l(D)$ or, equivalently,
\\[1ex]
$(a2)$
$(n - \lambda) L(\lambda , D^c) = (-1)^{n-1} \lambda L(n - \lambda , D)$.

Moreover, the matrices $L(D)$ and $L(D^c)$ are simultaneously diagonalizable.
\end{Theorem}

\bp~
Obviously, $D$ is balanced if and only if $D^c$ is balanced.
Let
$\{ \lambda _1,  \ldots, \lambda _n\}$ be the set of eigenvalues of $L(D)= L$.
Since $D$ is balanced, $L$ has an eigenvector
$X_n = J_{n1}$ with the eigenvalue $\lambda  _n = 0$, i.e. $LX_n = 0$.
Also since $D$ is balanced, by
Lemma  \ref{ortogLaplac},
there exists a basis $\{X_1,  \ldots, X_n\}$ of eigenvectors of $L$ such that $L X_i = \lambda_i X_i$ for every  $i = 1,\ldots, n$ and $X_i$ is orthogonal to $X_n = J_{n1}$ for every $i = 1,\ldots, n-1$.
Obviously, $L(K) = nI_n - J_{nn}$, where $K$ is the simple complete digraph with $n$ vertices.
Therefore $(nI_n - J_{nn} - L)X_n = 0$
and
 $ (nI_n - J_{nn} - L)X_i = (n -\lambda_i)X_i$ ~for~$i =1, 2, \ldots, n-1$.

Thus, $\lambda ^c_n = 0$ is an eigenvalue of $D^c$ corresponding to an eigenvector $X_n$ and $\lambda ^c_i =  n - \lambda_i $ is a Laplacian eigenvalue of $D^c$ corresponding to an eigenvector  $X_i$  for $i =1,  \ldots, n-1$.
\ep
\\[1ex]
\indent
The following Reciprocity Theorem is  true for all simple digraphs.
\begin{Theorem}
{\em \cite{KelLatvia,Ktree1}}
\label{di-reciprocityAllD}
Let  $D$ be a  simple  digraph with $n$ vertices.
Then
\\[1ex]
$(a1)$ there exists a bijection
$\sigma :
S'_l(D) \to S'_l(D^c)$ such that $x + \sigma ( x) = n$
for every $x \in S'_l(D)$ or, equivalently,
\\[1ex]
$(a2)$
$(n - \lambda) L(\lambda , D^c) = (-1)^{n-1} \lambda L(n - \lambda , D)$.
\end{Theorem}

\begin{example} Consider  a (simple) digraph $D$ with
\\
$V(D) = \{1, \ldots , 5\}$ and $E(D) = \{(1,2), (2,3), (3,4),
(4,5), (5,1),(1,3), (3,1)\}$.
\\
Then $S_l(D) = \{\lambda _0, \ldots , \lambda _4\}$ and
$S_l(D^c) = \{\lambda _1^c, \ldots , \lambda _5^c\}$, where $\lambda _5 = \lambda _5^c = 0$,

$\lambda _1 = 1.12256 - 0.744862~i $,
~$\lambda _1^c = 3.87744 + 0.744862~i $,

$\lambda _2 = 1.12256 + 0.744862~i $,
~$\lambda _2^c =  3.87744 - 0.744862~i $,

$\lambda _3 = 2.75488 $,
~$\lambda _3^c = 2.24512$,
~~and
$\lambda _4 = 2$,
~$\lambda _4^c = 3 $.
\\
Thus, $\lambda _s + \lambda _s^c = 5$ for every
$s \in \{1, 2, 3, 4\}$.
\end{example}
\indent
If $D$ is a digraph such that $d_{out}(x) = r$ for every $x \in V(D)$, then each
 $\alpha _i (D) = r - \lambda _i(D)$.
Therefore from Theorem \ref{di-reciprocityAllD}, we have   the following
\begin{corollary}
{\em \cite{Kcourse,KelLatvia}}
\label{reciprocity}
Let $D$ be a simple
digraph with $n$ vertices.
Suppose that $d_{out}(x) = r$ for every $x \in V(D)$.
Then
$$ A(\lambda, D^c) = (-1)^n
(\lambda - n + 1 + r)(\lambda +1 + r)^{-1}
A(-\lambda -1, D).$$
\end{corollary}

From Theorem \ref{di-reciprocityAllD} we also have a similar corollary  for digraphs with possible loops and no multiple arcs.
\begin{corollary}
{\em \cite{Kcourse,KelLatvia}}
\label{reciprocity'}
Let $D$ be an
$n$-vertex
digraph with possible loops.
Suppose that $d_{out}(x) = r$ for every $x \in V(D)$.
Then
$$A(\lambda, D_\circ^c) = (-1)^n
(\lambda - n + r)(\lambda + r)^{-1}
A(-\lambda, D).$$
\end{corollary}

\begin{lemma}\label{fAJ}
Let $D$ be an $r$-regular digraph with $n$ vertices
and let
$\{ \alpha _1,  \ldots, \alpha_n\}$ be the set of eigenvalues of $A(D)= A$, where $\alpha_n = r$.
 Let $f(x, y)$ be a polynomial with two variables and real coefficients.
Then matrix $f(A, J_{nn})$ has the eigenvalues $f(r, n)$ and $f(\alpha _i, 0 )$ for $i = 1,\cdots, n-1.$
\end{lemma}

\bp
Since $D$ is $r$-regular, $D$ is balanced and
$A(D) = rI - L(D)$. Therefore by
Lemma  \ref{ortogLaplac},
there exists a basis $\{X_1,  \ldots, X_n\}$ of eigenvectors of $A$ such that $A X_i = \alpha_i X_i$ for every  $i = 1,\ldots, n$ and $X_i$ is orthogonal to $X_n = J_{n1}$.
Then
\[f(A, J_{n n}) X_n = f(r, n) X_n
~~\mbox{\ and \ } ~~
f(A, J_{n n}) X_i= f(\alpha_i, 0) X_i~\mbox{\ for \ }
 i =1, \ldots , n-1.\]
\ep

\section{Adjacency spectra of   $D^{xyz}$ with $ z \in \{0, 1\} $}
\label{z-in-{0,1}}

\indent

Given a digraph $D$, we always denote by $A$, $T$, and $H$ the adjacency matrix,
the tail and the head incidence matrices of $D$, respectively. If $D$ is $r$-regular, then $e(D) = r v(D).$

\subsection{Spectra of   $D^{xyz}$ with $ z = 0$ and
$- \not \in \{x,y\}$}

\indent

Using Lemma \ref{PL} it is easy to prove the following theorem.
\begin{Theorem}
\label{xy0}
Suppose that one of the following holds:
\\[0.5ex]
$(c1)$  $D$ is a simple digraph and $x,y \in \{0,1,+, -\}$ or
\\[0.5ex]
$(c2)$ $D$ is a general digraph and
$x,y \in \{0,+\}$.

\noindent Then
$A(\lambda, D^{xy0}) =
A(\lambda, D^x) A(\lambda, (D^l)^y)$.
\end{Theorem}

Since $A(\lambda, D^+) = A(\lambda, D)$ and $A(\lambda, D^0) = \lambda^n$, we have from Theorem \ref{xy0} the following explicit formulas.
\begin{Theorem}
\label{xy0general}
Let  $D$ be a general digraph with $n$ vertices and $m$ arcs. Then
\\[0.5ex]
$(a1)$
$A(\lambda, D^{x00}) = \lambda^m A(\lambda, D^x)$ for $x\in \{0, +\}$,
\\[0.5ex]
$(a2)$
$A(\lambda , D^{0+0}) =  \lambda ^m  A(\lambda , D)$,  and
\\[0.5ex]
$(a3)$
$A(\lambda ,D^{++0}) =
\lambda^{m - n} A(\lambda , D)^2$.
\end{Theorem}

From $(a1)$ and $(a2)$ in Theorem \ref{xy0general} we have a simple but very interesting observation.
\begin{corollary}
\label{Sp(+00=0+0)}
Let $D$ be a general digraph. Then
$D^{+00} \sim ^A D^{0+0}$ and if $D \not \cong D^l$, then $D^{+00} \not \cong D^{0+0}$.
\end{corollary}

If $D$ is a simple digraph, then
$A(\lambda, D^{1})= (\lambda + 1)^{n-1} (\lambda - n + 1)$. Therefore from Lemma \ref{PL} we have:
\begin{Theorem}
\label{xy0-no-r}
Let $D$ be a simple digraph with $n$ vertices and $m$ arcs.
  Then
\\[0.5ex]
$(a1)$
$ A(\lambda, D^{x10}) =
 (\lambda + 1)^{m-1} (\lambda - m + 1) A(\lambda, D^x)$ for $x \in \{0, 1, +\}$,
\\[0.5ex]
$(a2)$
$A(\lambda, D^{100}) = \lambda^{m} (\lambda + 1)^{n-1} (\lambda - n +1)$, and
\\[0.5ex]
$(a3)$
$A(\lambda ,D^{1+0}) =
\lambda^{m - n} (\lambda + 1)^{n-1} (\lambda - n + 1)A(\lambda , D)$.
\end{Theorem}

\subsection{Spectra of  $D^{xyz}$ with $ z = 0$ and
$- \in \{x,y\}$
 for a simple regular digraph $D$}

\indent

From Lemma \ref{PL}, Corollary \ref{reciprocity}, and
Theorems
\ref{xy0general} and \ref{xy0-no-r} we obtain the adjacency polynomials of the other digraphs $D^{xyz}$ with $z=0$.

\begin{Theorem}
\label{xy0+r}
Let $D$ be a simple $r$-regular
digraph with $n$ vertices and $m$ arcs.
  Then
\\[1.5ex]
$(a1)$
$A(\lambda, D^{-00}) = \lambda^m A(\lambda, D^c) = (-1)^n
\lambda^{m} (\lambda - n + r + 1)(\lambda + r + 1)^{-1}
A(-\lambda -1, D)$,
\\[1.5ex]
$(a2)$
 $A(\lambda, D^{-10}) = (-1)^n (\lambda +1)^{m-1} (\lambda- m + 1)
  (\lambda - n + r + 1)(\lambda + r + 1)^{-1}
 A(-\lambda -1, D)$,
\\[1.5ex]
$(a3)$
$A(\lambda, D^{0-0}) = (-1)^n  \lambda^{n} (\lambda + 1)^{m-n}
(\lambda - m + r + 1)(\lambda + r + 1)^{-1}
A(-\lambda - 1, D)$,
\\[1.5ex]
$(a4)$
$A(\lambda, D^{1-0}) = (-1)^n (\lambda + 1)^{m - 1} (\lambda -n + 1)
(\lambda - m + r + 1)(\lambda + r + 1)^{-1}
 A(-\lambda -1, D)$,
\\[1.5ex]
$(a5)$
 $A(\lambda ,D^{-+0}) =
(-1)^n \lambda^{m - n}
(\lambda - n + r + 1)(\lambda + r + 1)^{-1}
A(- \lambda - 1, D) ~A(\lambda , D) $,
\\[1.5ex]
$(a6)$
$A(\lambda , D^{+-0}) =  (-1)^n
 (\lambda + 1)^{m - n}
(\lambda - m + r + 1)(\lambda + r + 1)^{-1}
A(- \lambda - 1, D) ~A( \lambda , D)$, and
\\[1.5ex]
$(a7)$
 $A(\lambda ,D^{--0}) =
( \lambda + 1)^{m - n}
(\lambda - n + r + 1)(\lambda - m + r + 1)
(\lambda + r + 1)^{-2}
A( - \lambda - 1 , D)^2$.
\end{Theorem}

\subsection{Spectra of   $D^{xyz}$ with
$z = 1$ and $x,y \in \{0, 1\}$ for a simple digraph $D$}

\begin{Theorem}
\label{simple,z1}
Let $D$ be a simple digraph with $n$ vertices and $m$ arcs.
  Then
\\[1.5ex]
$(a1)$  $A(\lambda, D^{001}) = \lambda^{m + n - 2} (\lambda^2 - mn)$,
\\[1.5ex]
$(a2)$  $A(\lambda , D^{101}) =  \lambda^{m-1} (\lambda + 1)^{n - 1}
    (\lambda^2 + \lambda - n \lambda - mn), $
\\[1.5ex]
$(a3)$  $A(\lambda, D^{011}) = \lambda^{n-1} (\lambda + 1)^{m-1} (\lambda^2 + \lambda - m\lambda - mn)$, and
\\[1.5ex]
$(a4)$  $A(\lambda, D^{111}) = (\lambda + 1)^{m+ n - 1} (\lambda - m - n + 1)$.
\end{Theorem}

\bp~
We give only the proof of claim $(a2)$.
The proofs of the remaining claims
 are similar.
The adjacency matrix of $D^{101}$ is

\[ A(D^{101})=
\left(
     \begin{array}{cc}
       J_{nn} - I_n  & J_{nm} \\
       J_{mn} & 0 \\
     \end{array}
   \right).
\]
Therefore
\[ A(\lambda, D^{101}) = \left|
     \begin{array}{cc}
       (\lambda + 1) I_n - J_{nn} & -J_{nm} \\
       -J_{mn} & \lambda I_m \\
     \end{array}
   \right|. \]
Using Lemmas \ref{lemABCD} and \ref{fAJ} and assuming that $\lambda \neq 0$, we  obtain:
\begin{eqnarray*}
 A(\lambda, D^{101}) &=& \lambda^{m-n} |(\lambda^2 + \lambda)I_n - J_{nn}(\lambda + m)| \\
   &=& \lambda^{m-n}(\lambda^2 + \lambda - n(\lambda + m))
   \prod_{i=1}^{n-1}\{\lambda^2 + \lambda\} \\
   &=& \lambda^{m-1} (\lambda + 1)^{n-1}(\lambda^2 + \lambda - n \lambda - mn).
\end{eqnarray*}
\ep

 \subsection{Spectra of   $D^{xyz}$ with $ z = 1$,
 $xy \in \{+0, 0+, ++\}$ for a general regular digraph $D$}

\begin{Theorem}
\label{+01}
Let $D$ be a general  $r$-regular digraph with $n$ vertices and $m$ arcs.
  Then
  \[A(\lambda , D^{+01}) =  \lambda^{m-1}
    (\lambda^2 -  r\lambda - m n )(\lambda - r)^{-1}
  A(\lambda, D). \]
\end{Theorem}

\bp~
Let $S_a(D) = \{\alpha _i: i = 1, \ldots, n\}$, where
$\alpha _n = r$.
The adjacency matrix of $D^{+01}$ is
\[ A(D^{+01}) =\left(
     \begin{array}{cc}
        A & J_{nm} \\
       J_{mn} & 0 \\
     \end{array}
   \right).
\]
Then
\[ A(\lambda, D^{+01}) = \left|
     \begin{array}{cc}
       \lambda I_n - A & -J_{nm} \\
       -J_{mn} & \lambda I_m \\
     \end{array}
   \right|. \]
 Using Lemmas \ref{lemABCD} and \ref{fAJ} and assuming that $\lambda \neq 0$, we obtain:
\begin{eqnarray*}
  A(\lambda, D^{+01}) &=& \lambda^{m-n} |\lambda^2 I_n - \lambda A - m J_{nn}| \\
   &=& \lambda^{m-n}(\lambda^2 -  r\lambda - m n)
   \prod_{i=1}^{n-1}\{\lambda^2 - \alpha_i  \lambda\}.
\end{eqnarray*}
Therefore
\[
A(\lambda , D^{+01}) = \lambda^{m-1}
    (\lambda^2 -  r\lambda - m n )(\lambda - r)^{-1}
  A(\lambda, D).
\]
\ep

\begin{Theorem}\label{0+1}
Let $D$ be
a general  $r$-regular digraph with $n$ vertices and $m$ arcs.
  Then
  \[A(\lambda , D^{0+1}) =  \lambda^{m-1}
  (\lambda^2 -  r\lambda - m n )(\lambda - r)^{-1}
  A(\lambda, D). \]
\end{Theorem}

\bp~
Let $S_a(D^l) = \{\alpha _i: i = 1, \ldots, m\}$, where
$\alpha _m = r$.
The adjacency matrix of $D^{0+1}$ is
\[A(D^{0+1}) = \left(
     \begin{array}{cc}
        0 & J_{nm} \\
       J_{mn} & A(D^l) \\
     \end{array}
   \right).
\] Then we have:
\[ A(\lambda, D^{0+1}) = \left|
     \begin{array}{cc}
       \lambda I_n & -J_{nm} \\
       -J_{mn} & \lambda I_m - A(D^l) \\
     \end{array}
   \right|. \]
 Using Lemmas \ref{lemABCD} and \ref{fAJ} and assuming that $\lambda \neq 0$, we obtain:
\begin{eqnarray*}
 A(\lambda, D^{0+1}) &=& \lambda^{n - m} |\lambda^2 I_m - \lambda A(D^l) - n J_{mm}| \\
   &=& \lambda^{n-m}(\lambda^2 -  r\lambda - m n)
   \prod_{i=1}^{m-1}\{\lambda^2 - \alpha_i  \lambda\}.
\end{eqnarray*}
Hence
   \[ A(\lambda, D^{0+1}) =  \lambda^n
  (\lambda^2 - r \lambda - m n )(\lambda^2 - r \lambda)^{-1}
 A(\lambda, D^l). \]
Now by Lemma \ref{PL},
\[  A(\lambda, D^{0+1}) =  \lambda^{m-1}
(\lambda^2 -  r\lambda - m n )(\lambda - r)^{-1}
 A(\lambda, D).  \]
\ep

\begin{figure}[ht]
\begin{center}
\scalebox{0.7}[.7]{\includegraphics{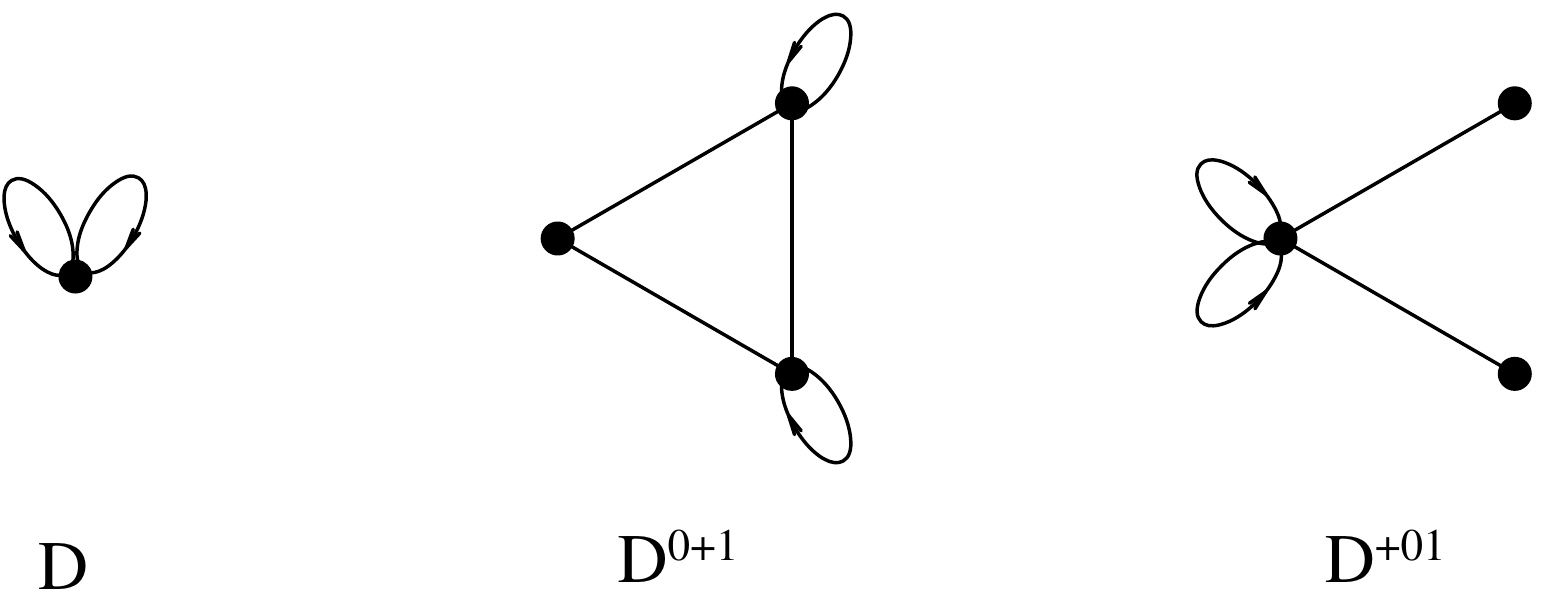}}
\end{center}
\caption{Cospectral and non-isomorphic digraphs $D^{0+1}$ and $D^{+01}$, where each undirected edge represents a pair of arcs with opposite directions.}\label{0+1=+01}
\end{figure}

From Theorems \ref{+01} and \ref{0+1} we have:
\begin{corollary}
\label{Sp(+01=0+1)}
If $D$ is a general $r$-regular digraph with $n$ vertices and $m$ arcs, then $D^{+01}$ and $D^{0+1}$ are cospectral and they have $m-1$ eigenvalues equal to zero, $n -1$ eigenvalues $\alpha_i$,
$i = 1, \ldots , n-1$, and two  additional eigenvalues
$ \frac{1}{2} (r \pm \sqrt{r^2 + 4 m n}).$ Moreover, if $r > 1$, then $D^{+01} \not \cong D^{0+1}$.
\end{corollary}

\bp~
We prove the last claim. If $r > 1$, we have $m > n$. Then the maximum size of independent sets in $D^{+01}$ is $m$ while this number in $D^{0+1}$ is $n$. Therefore  $D^{+01} \not \cong D^{0+1}$.
\ep

For example, let $D$ be the digraph with one vertex and $n-1$ loops. Then $D^{0+1}$ and $D^{+01}$ are cospectral but their underlying graphs are not isomorphic  (see Figure \ref{0+1=+01} for $n = 3$).
\\[1ex]
\indent
The proof of the next theorem is similar to those of Theorems  \ref{+01} or \ref{0+1}.
\begin{Theorem}\label{++1}
Let $D$ be
a general $r$-regular digraph with $n$ vertices and $m$ arcs.
  Then
  \[A(\lambda , D^{++1}) =  \lambda^{m-n}
  ((\lambda - r)^2 - m n )(\lambda - r)^{-2}
  A(\lambda, D)^2. \]
\end{Theorem}

\subsection{Spectra of   $D^{xyz}$ with
$z = 1$ and $\{x,y \}= \{1,+\}$ or $- \in \{x,y\}$ for a simple regular digraph $D$}
\begin{Theorem}
\label{1+1,+11}
Let $D$ be a simple $r$-regular digraph with $n$ vertices and $m$ arcs.
  Then
\\[1.5ex]
$(a1)$  $A(\lambda , D^{1+1}) =  \lambda^{m-n} (\lambda + 1)^{n-1}
    ((\lambda - r)(\lambda + 1 - n) - mn))(\lambda - r)^{-1}
  A(\lambda, D),$ and
  \\[1.5ex]
$(a1)$ $A(\lambda , D^{+11}) =  (\lambda + 1)^{m-1}
    ((\lambda - r)(\lambda + 1 - m) - mn)(\lambda - r)^{-1}
  A(\lambda, D).$
\end{Theorem}

The proof of the above theorem is similar to those for Theorems \ref{simple,z1} or \ref{+01}.

\begin{Theorem}
\label{-01}
Let $D$ be a simple $r$-regular digraph with $n$ vertices and $m$ arcs.
  Then
  \[A(\lambda , D^{-01}) = (-1)^n \lambda^{m-1}
    (\lambda (\lambda + 1 + r - n) - mn)(\lambda + 1 +  r)^{-1}
  A(-\lambda - 1, D). \]
\end{Theorem}

\bp~
Let $S_a(D) = \{\alpha _i: i = 1, \ldots, n\}$ and
$\alpha _n = r$.
The adjacency matrix of $D^{-01}$ is

\[ A(D^{-01})~=~
\left(
     \begin{array}{cc}
       J_{nn} - I_n - A &~~ J_{nm} \\
       J_{mn} &~~ 0 \\
     \end{array}
   \right).
\]
Then
\[ A(\lambda, D^{-01}) ~=~ \left|
     \begin{array}{cc}
       (\lambda + 1) I_n - J_{nn} + A & ~~-J_{nm} \\
       -J_{mn} & ~~~~\lambda I_m \\
     \end{array}
   \right|. \]
Using Lemmas \ref{lemABCD} and \ref{fAJ} and assuming that $\lambda \neq 0$, we obtain:
\begin{eqnarray*}
 A(\lambda, D^{-01}) &=& \lambda^{m-n} |(\lambda^2 + \lambda)I_n + \lambda A + J_{nn}(-\lambda - m)| \\
   &=& \lambda^{m-n}(\lambda^2 + \lambda + \lambda r - n(\lambda + m))
   \prod_{i=1}^{n-1}\{\lambda^2 + \lambda  + \lambda ~\alpha_i\}.
\end{eqnarray*}
Therefore
\begin{eqnarray*}
 A(\lambda , D^{-01}) &=& (-1)^n \lambda^m
    (\lambda^2 + \lambda + \lambda r - n (\lambda + m))(\lambda^2 + \lambda + \lambda r)^{-1}
  A(-\lambda - 1, D) \\
&=& (-1)^n \lambda^{m-1}
    (\lambda (\lambda + 1 + r - n) - mn) (\lambda + 1 +  r)^{-1}
  A(-\lambda - 1, D).
\end{eqnarray*}
\ep

\begin{corollary}
\label{Sp(-01)}
Let $D$ be a simple $r$-regular digraph with $n$ vertices and $m$ arcs and
\\
$S_a(D) = \{\alpha _i: i = 1, \ldots, n\}$, where $\alpha _n = r$.
Then $D^{-01}$ has $m-1$ eigenvalues equal to zero,  $n-1$  eigenvalues  $- (\alpha_i +1)$, where $i = 1, \ldots , n-1$, and two additional eigenvalues
\[ \frac{1}{2}( n - r - 1 \pm \sqrt{(n - r - 1)^2 + 4 r n^2 }). \]
\end{corollary}

The adjacency polynomials of
$D^{xyz}$ for the other cases when
$z = 1$ and $- \in \{x,y\}$ can be found in Appendix.
The proofs of these cases are using Lemmas \ref{PL}, \ref{AL}, \ref{lemABCD}, and \ref{fAJ} and are similar to that of  Theorem \ref{-01}.

\section{Adjacency spectra of   $D^{xyz}$ with
$z \in \{ +, -\}$}
\label{z-in-{+,-}}
\indent

In this section we consider mainly the adjacency spectra of $D^{xyz}$ for $z\in \{ +, - \}$ and $|\{x, y\} \cap \{+, -\}| = 1$. We also describe the spectrum of $D^{++-}$ for a general regular digraph $D$. The spectra of $D^{+++}, D^{00+}, D^{+0+}$ and $D^{0++}$ for a general digraph $D$ can be found in \cite{zhang}. The spectra of the other $D^{xyz}$ with $x, y, z \in \{+, -\}$ for simple regular digraph $D$ are given
in \cite{liumeng}.
All these formulas are also included in Appendix.

\subsection{Spectra of   $D^{xyz}$ with
$z = +$ and  $| \{x, y \} \cap \{0, +\}| \leq 1$ for a simple regular digraph $D$}

\begin{Theorem}\label{11+}
Let $D$ be
a simple $r$-regular digraph with $n$ vertices and
$m$ arcs.
Then
\[ A( \lambda , D^{11+}) = (\lambda + 1)^{m-n}~
\frac{(\lambda + 1)^2 - r - n((\lambda + 1)(r + 1) - m)}{(\lambda + 1)^2 - r}~
A((\lambda + 1)^2, D).
\]
\end{Theorem}

\bp~
Let $S_a(D) = \{\alpha _1, \ldots , \alpha _n\}$, where
$\alpha _n = r$.
The adjacency matrix of $D^{11+}$ is
\[A(D^{11+})= \left( \begin{array}{cc}  J_{nn}-I_n     &~~ T
\\[1ex]
H^{\top}  &~~ \quad J_{mm} - I_m
\end{array}
\right).\]
Then
\begin{eqnarray*}
 A( \lambda , D^{11+}) &=& \left| \begin{array}{cc}  (\lambda + 1) I_n - J_{nn}     &~~ - T
\\[1ex]
-H^{\top}  &~~ \quad
(\lambda + 1) I_m - J_{mm}
\end{array}
\right|.
\end{eqnarray*}
By Lemma \ref{lemJ} $(a2)$, $J_{mn}T = J_{mm}$. Hence multiplying the first row of the above block determinant by $-J_{mn}$ and adding the result to the second row, we obtain:
\[
A( \lambda , D^{11+}) =  \left| \begin{array}{cc}  (\lambda + 1) I_n - J_{nn}     &~~ - T
\\[1ex]
-H^{\top} -J_{mn}(\lambda + 1 - n) &~~ \quad
(\lambda + 1) I_m
\end{array}
\right|. \]
Now assuming that $\lambda  \neq -1$ and using Lemmas \ref{AL} $(a1)$, \ref{lemJ} $(a1)$, \ref{lemABCD}, and \ref{fAJ},   we obtain:
\begin{eqnarray*}
 A( \lambda , D^{11+})
&=& (\lambda + 1)^{m-n} |(\lambda + 1)((\lambda + 1) I_n - J_{nn}) - T(H^{\top} + J_{mn}(\lambda + 1 - n))| \\
&=& (\lambda + 1)^{m-n} |(\lambda + 1)^2 I_n - A - J_{nn}(\lambda + 1 + r(\lambda + 1 - n))| \\
&=& (\lambda + 1)^{m-n} ((\lambda + 1)^2 - r - n((\lambda + 1)(r+ 1) - rn)) \prod_{i=1}^{n-1} \{(\lambda + 1)^2 - \alpha_i\}.
\end{eqnarray*}
Note that $m = nr$. Hence
\[ A( \lambda , D^{11+}) =
 (\lambda + 1)^{m - n}~
 \frac{(\lambda + 1)^2 - r - n((\lambda + 1)(r + 1) - m)}{(\lambda + 1)^2 - r}~
  A((\lambda + 1)^2, D). \]
\ep

The remaining situations when $z = +$ and  $| \{x, y \} \cap \{0, +\}| \leq 1$ can be considered similarly
(see Appendix).

\subsection{Spectra of   $D^{xyz}$ with
$z = -$, $x,y \in \{0,+\}$ for a general regular digraph $D$}

\begin{Theorem}
\label{++-}
Let $D$ be  a general $r$-regular digraph with $n$ vertices and $m$ arcs. Then
\begin{eqnarray*}
 A( \lambda ,D^{++-})  &=& \lambda^{m - n}
(\lambda^2 - r(2 \lambda + 1) + r^2 + 2m - mn)(\lambda^2 - r(2 \lambda + 1) + r^2)^{-1}
\\[0.5ex]
   & &  A(
  2^{-1}
   (2\lambda + 1 + \sqrt{4\lambda + 1}) , D)  ~A(
 2^{-1}   (2\lambda + 1 - \sqrt{4\lambda + 1}) , D).
\end{eqnarray*}

\end{Theorem}

\bp
Let $S_a(D) = \{\alpha _1, \ldots , \alpha _n\}$, where
$\alpha _n = r$.
The adjacency polynomial  of $D^{++-}$ is
\[A( \lambda, D^{++-})= \left| \begin{array}{cc} \lambda I_n - A
&~~ - J_{nm} + T
\\[1ex]
- J_{mn} + H^{\top}  &~~ \quad
\lambda I_m - A^l
\end{array}
\right|.\]
Multiplying the first row of the above block matrix by
$(H^\top + (1-n)^{-1}J_{mn})$
and adding the result to the second row, we obtain:

\[A( \lambda, D^{++-})= \left| \begin{array}{cc} \lambda I_n - A
&~~ - J_{n m} + T
\\[1ex]
 H^{\top}((\lambda + 1)I_n - A) - J_{mn} +
((1-n)^{-1}
 J_{mn})(\lambda I_n - A)  &~~ \quad
\lambda I_m
\end{array}
\right|.\]
 Using Lemmas \ref{AL}, \ref{lemJ} $(a5)$, and \ref{lemABCD} and assuming that $\lambda \ne  0$, we have:
\begin{eqnarray*}
A( \lambda , D^{++-}) & = &
\lambda^{m-n} |\lambda^2 I_n - \lambda A + (J_{nm} - T)( H^\top((\lambda + 1) I_n - A) + J_{mn}( -1 + \frac{\lambda - r}{1 - n}))|
\end{eqnarray*}
Again using Lemmas \ref{AL}, claims $(a1)$, $(a3)$, and $(a5)$ of \ref{lemJ},  and   \ref{fAJ}, we obtain noting that $m = nr$:
\begin{eqnarray*}
A( \lambda , D^{++-}) & = &
 \lambda^{m-n} |\lambda^2 I_n - (2\lambda + 1) A + A^2 + (2-n) r J_{nn}|
\\[1ex]
& = & \lambda^{m - n} \frac{\lambda^2 - (2\lambda + 1)r + r^2 +2m - m n}{\lambda^2 - (2\lambda + 1)r + r^2}
\prod_{i=1}^{n} \{\lambda^2 - (2\lambda + 1) \alpha_i + \alpha_i^2\}
\\[1ex]
& = & \lambda^{m - n} \frac{\lambda^2 - (2\lambda + 1)r + r^2 +2m - m n}{\lambda^2 - (2\lambda + 1)r + r^2} \times
\\[1.5ex] & &
 A(
2^{-1} (2\lambda + 1 + \sqrt{4 \lambda + 1}), D)
~A(
2^{-1}(2\lambda + 1 - \sqrt{4 \lambda + 1}), D).
\end{eqnarray*}
\ep

Similarly, we can obtain the following result for general regular digraphs.
\begin{Theorem}
\label{{0,+}z-}
Let $D$ be a general $r$-regular digraph with $n$ vertices and $m$ arcs. Then
\\[1.5ex]
$(a1)$
$ A( \lambda, D^{00-}) =
 \lambda^{m-n}
  (\lambda^2 - r + 2m - mn)(\lambda^2 - r)^{-1}
 A(\lambda^2, D)  $ and
\\[1.5ex]
$(a2)$
$ A( \lambda, D^{+0-}) = A( \lambda, D^{0+-})$

 $= \lambda^{m - n} (\lambda + 1)^{n}
   (\lambda^2 - r(\lambda + 1) + 2m - mn)
  (\lambda^2 - r (\lambda + 1))^{-1}
 A(\frac{\lambda^2}{\lambda + 1}, D).  $
\end{Theorem}

From $(a2)$ in  Theorem \ref{{0,+}z-} we have the following result  similar to Corollaries \ref{Sp(+00=0+0)} and
\ref{Sp(+01=0+1)}.
\begin{corollary}
\label{Sp(+0-=0+-)}
Let $D$ be a general $r$-regular digraph.
Then $D^{+0-}$ and $D^{0+-}$ are cospectral digraphs and $D^{+0-} \not \cong D^{0+-}$ for $r>1$.
\end{corollary}

\subsection{Spectra of   $D^{xyz}$ with
$z = -$, $|\{x,y\} \cap \{0,+\}| \leq 1$ for a simple regular digraph $D$}

\begin{Theorem}
\label{-0-}
Let $D$ be a simple $r$-regular digraph with $n$ vertices and $m$ arcs. Then
 \[ A( \lambda, D^{-0-}) = \lambda^{m-n}(1-\lambda)^n~
\frac{\lambda^2 - \lambda (n - r - 1) + 2m - mn - r}
 {\lambda^2 + \lambda (r+1)-  r}~
A(\frac{\lambda^2 + \lambda}{1 - \lambda}, D).  \]
\end{Theorem}

\bp
Let $S_a(D) = \{\alpha _1, \ldots , \alpha _n\}$, where
$\alpha _n = r$.
The adjacency matrix of $D^{-0-}$ is
\[A(D^{-0-}) = \left(
    \begin{array}{cc}
      J_{nn} - I_n - A &~~  J_{nm} - T \\[1ex]
      J_{mn} - H^{\top} &~~ 0 \\
    \end{array}
  \right).
\]
Thus, \[A(\lambda, D^{-0-}) = \left|
    \begin{array}{cc}
     (\lambda + 1) I_n - J_{nn} + A &~~ - J_{nm} + T \\[1ex]
      -J_{mn} + H^{\top} &~~ \lambda I_m   \\
    \end{array}
  \right|.\]
Using Lemmas \ref{AL} $(a1)$, \ref{lemJ} $(a1)$, $(a3)$, and \ref{lemABCD},  and assuming that $\lambda \neq 0$, we obtain:
 \[A(\lambda, D^{-0-}) = \lambda^{m-n} |(\lambda^2 + \lambda)I_n - \lambda J_{nn} + \lambda A - (m - 2r) J_{nn} - A|. \]
 Now using Lemma \ref{fAJ} we obtain:
 \[A(\lambda, D^{-0-}) = \lambda^{m-n}(\lambda^2 + \lambda + r(\lambda - 1) - n (\lambda + m - 2r))
 \prod_{i=1}^{n-1} (\lambda^2 + \lambda + \alpha_i(\lambda - 1)).\]
Therefore
\[A(\lambda, D^{-0-}) = \lambda^{m-n}(1-\lambda)^n~
\frac{\lambda^2 - \lambda (n - r - 1) + 2m - mn - r}
{\lambda^2 + \lambda (r+1)-  r}~
A(\frac{\lambda^2 + \lambda}{1 - \lambda}, D).\]
\ep
\\[1ex]
\indent
The formulas for
the remaining cases when
$z = -$ and $|\{x,y\} \cap \{0,+\}| \leq 1$
can be found similarly (see Appendix).

\section{Digraph-functions and their $xyz$-transformations}

\label{digraph-functions}
\indent

A digraph $D = (V, E)$ is called a {\em digraph-function} \cite{Kcourse}
if  there exists a function
\\
$f: V \to V$ such that
$(x,y) \in E$ if and only if $y = f(x)$. Similar digraphs were considered in \cite{H&S96}.

A digraph $D$ is called a  {\em directed cycle} or simply, {\em dicycle} if $D$ is  connected and
$d_{in}(x) = d_{out}(x) = 1$ for every $x \in V(D)$.
A {\em directed $xy$-path } or simply, {\em $xy$-dipath} is a  digraph obtained from a dicycle with an arc $(y,x)$ by removing  arc $(y,x)$.

Let ${\cal F}$ denote the set of digraphs $F$ such that
each component of $F$ is either a digraph-function or its inverse. Let ${\cal CF}$ denote the set of connected digraphs in ${\cal F}$.
It is easy to see that if $F \in {\cal CF}$,
then $F$ has a unique directed cycle $C$
(possibly, a loop). Let $c(F) = v(C)$.

\subsection{Digraph-function criteria}
\label{digraph-function-criteria}

\indent

The following simple observation provides different digraph-function criteria.
\begin{lemma}
\label{crt}
Let $D$ be a connected digraph with $n$ vertices. Then the following statements are equivalent:
\\[1ex]
$(a1)$ $D$ is a digraph-function $($resp., the inverse of digraph-function$)$,
\\[1ex]
$(a2)$ each vertex of $D$ has out-degree one $($resp., in-degree one$)$,
\\[1ex]
$(a3)$ $t: E(D) \to V(D)$ $($resp., $h: E(D) \to V(D)$$)$ is an isomorphism from $D^l$ to $D$,
\\[1ex]
$(a4)$
if $V(D) = \{1, \ldots , n\}$ and $e_i = t ^{-1}(i)$, then $T(D) = I_n$ and
$H^\top (D) = A(D)$
$($resp., $H(D) = I_n$  and $T(D) = A(D)$, and
$
)$, and $A(D^l) = H^\top (D) = A(D)$,
\end{lemma}

Here is another interesting criterion for a digraph to be a digraph-function or its inverse.
\begin{Theorem}\label{Dl=D}
Let $D$ be a digraph. Then  $D \in {\cal F}$
 if and only if  $D$ is isomorphic to $D^l$.
\end{Theorem}

\bp
By $(a3)$ in Proposition \ref{crt},
 if $D \in {\cal F}$, then  $D$ is isomorphic to
$D^l$. We will prove that if $D$ is isomorphic to
$D^l$, then $D \in {\cal F}$.
It is sufficient to prove our claim for a connected digraph $D$.
By $(a2)$ in Proposition \ref{crt}, it is sufficient to prove, that each vertex of $D$ has out-degree one or each vertex of $D$ has in-degree one. 

Suppose, on the contrary, that $D$ has a vertex $u$ such that $d_{out}(u) \ge 2$ or $d_{in}(u) \ge 2$. We can assume that $d_{out}(u) \ge 2$ (the case $d_{in}(u) \ge 2$ can be considered similarly).

First we define special digraphs which we call {\em claws}.
Given three disjoint dipaths $p'Pp$, $q'Qq$,  and $rRr'$ with $p \ne p'$ and $q \ne q'$, let
digraph $Y$ be obtained from $p'Pp$, $q'Qq$,  and $rRr'$ by identifying three vertices $p'$, $q'$, and $r'$ with a new vertex $c$. We call
$Y$ an {\em $(p,q,r)$-claw with the center} $c$ or simply, a {\em claw}. If, in particular, $r = r'$, we call $Y$ a {\em 2-leg claw}.
Let ${\cal Y}$ and ${\cal Y}^l$ be the sets of all claw sub-digraphs in $D$ and $D^l$, respectively.

Since $d_{out}(u) \ge 2$, clearly  $D$ has a 2-leg claw $T$ with center $u$ as a  sub-digraph.
Therefore   ${\cal Y} \ne \emptyset $. Since $D$ is isomorphic to $D^l$, also ${\cal Y}^l \ne \emptyset $. Given a sub-digraph 
$F$ of $D^l$, let $F^{-l}$ denote the sub-digraph of $D$ such that 
$(F^{-l})^l = F$.

Let $S$ and $Z$ be largest  claws (i.e. claws with the maximum number of arcs) in $D$ and $D^l$, respectively. Let $m(D) = e(S)$ and  $m(D^l) = e(Z)$. %
Since $D$ and $D^l$ are isomorphic, clearly $m(D) = m(D^l)$.
Obviously, $e(Z^{-l}) = e(Z) +1$ and $Z^{-l} \in {\cal Y}$. 
Thus, $m(D) \ge e(Z^{-l}) = e(Z) +1 > e(Z) = m(D^l)$, a contradiction.
\ep

\subsection{Spectra of digraph-function $xyz$-transformations for $z = 0$ or
$z \in \{1, +\}$    and $x,y \in \{0, z\}$}
\label{df-spectra}

\indent

In this subsection we will
 consider the triples $xyz$ such that $z = 0$
or $z \in \{1, +\}$  and  $x, y \in \{0, z\}$, and describe $A(\lambda ,F^{xyz})$ for  every connected graph-function $F$ and its inverse in terms of the spectrum of
$F$ (i.e., in terms of $v(F)$ and $c(F))$.

\begin{Theorem}
\label{A(x,F)}
Let
$F \in \cal{CF}$,
$C$ be the directed cycle in $F$,
$v(F) = n$, and $c(F) = v(C) = k$. Then
 \\[0.7ex]
 $(a1)$
$A(\lambda , F) = \lambda ^{n - k} A(\lambda , C)
=  \lambda ^{n - k} (\lambda ^k - 1)$ and
 \\[0.7ex]
 $(a2)$
 $A(\lambda, F^c) = (\lambda + 1)^{n - k} (\lambda - n + 2) (\lambda + 2)^{-1} ((\lambda + 1)^k - (-1)^k).$
\end{Theorem}

\bp
Claim $(a1)$ is obvious. We prove
$(a2)$.
By Lemma \ref{Inverse} $(a1)$, it is sufficient to prove our claim for a
digraph-function $F$. Then $d_{out}(x,F) =1$ for every $x \in V(F)$. Now  $(a2)$ follows from $(a1)$ by Corollary \ref{reciprocity} with $r = 1$.
\ep

\vskip 1ex
Now we can easily describe $A(\lambda, F^{xyz})$ for $F  \in \cal{CF}$ in terms of
 $v(F)$ and $c(F)$ when  $z=0$.
\begin{Theorem}
\label{Fz=0}
Let
$F  \in \cal{CF}$, $v(F) = n$, and $c(F) =  k$.
 Then
 \\[0.7ex]
 $(a1)$
$A(\lambda , F^{000}) =  \lambda ^{2n},$
 \\[0.7ex]
 $(a2)$
$A(\lambda , F^{100}) = A(\lambda, F^{010}) = (\lambda - n + 1) \lambda ^n (\lambda + 1)^{n - 1}$,
 \\[0.7ex]
 $(a3)$
$A(\lambda , F^{110}) = (A(\lambda , F^1))^2 = (\lambda - n + 1)^2 (\lambda + 1)^{2n - 2}$,
 \\[0.7ex]
 $(a4)$
$A(\lambda , F^{+00}) = A(\lambda , F^{0+0}) = A(\lambda , F^0)~A(\lambda , F)=
\lambda ^n A(\lambda , F) =
\lambda ^{2n - k} (\lambda ^k - 1)$,
 \\[0.7ex]
 $(a5)$
$A(\lambda , F^{++0}) = (A(\lambda , F))^2 =
 \lambda ^{2(n - k)} (\lambda ^k - 1)^2$,
 \\[0.7ex]
$(a6)$
$A(\lambda , F^{+10}) = A(\lambda , F^{1+0}) = A(\lambda , F^1)~A(\lambda , F) =
(\lambda - n + 1) \lambda^{n - k} (\lambda + 1)^{n - 1} (\lambda^k - 1)$,
 \\[0.7ex]
$(a7)$
$A(\lambda , F^{-00}) = A(\lambda , F^{0-0}) =  \lambda^n (\lambda + 1)^{n-k} (\lambda - n + 2)
  (\lambda + 2)^{-1} ((\lambda + 1)^k - (-1)^k )$,
 \\[0.7ex]
$(a8)$
$A(\lambda , F^{--0}) = (A(\lambda , F^c))^2 = (\lambda + 1)^{2n - 2k} (\lambda - n + 2)^2
  (\lambda + 2)^{-2} ((\lambda + 1)^k - (-1)^k)^2$,
 \\[0.7ex]
$(a9)$
$A(\lambda , F^{-10}) = A(\lambda , F^{1-0}) = A(\lambda , F^1)~A(\lambda , F^c)$
\\
\indent ~ $ = (\lambda - n + 1) (\lambda - n + 2) (\lambda + 2)^{-1} (\lambda + 1)^{2n - k - 1}
 ((\lambda + 1)^k - (-1)^k )$, and
 \\[0.7ex]
$(a10)$
$A(\lambda , F^{+-0}) =
A(\lambda , F^{-+0}) = A(\lambda , F) ~A(\lambda , F)^c $
\\
\indent ~ $ = (\lambda - n + 2) (\lambda + 2)^{-1}  \lambda^{n - k}
(\lambda^k - 1) (\lambda + 1)^{n - k}  ((\lambda + 1)^k - (-1)^k )$.
\end{Theorem}

\bp
Obviously, $F^{++0} = F\cup F^l$.
By Theorem \ref{crt}, $F$ and $F^l$ are isomorphic. Therefore it is easy to see that claims $(a1)$ - $(a3)$ are true.
The other cases follow directly from Theorem \ref{A(x,F)}.
\ep
\\[1ex]
\indent
It is also easy to describe $A(\lambda , F^{xyz})$ for $F \in \cal{F}$
when $z = 1$ and $x, y \in \{0, 1\}$. Actually, it is just a special case of Theorem \ref{simple,z1}.
\begin{Theorem}
 \label{Fz1xy01}
Let
$F  \in \cal{F}$ and $v(F) = n$.
 Then
  \\[0.7ex]
$(a1)$
$A(\lambda , F^{001}) = \lambda ^{2n - 2} (\lambda ^{2} - n^2)$.
 \\[0.7ex]
$(a2)$
$A(\lambda , F^{011}) = A(\lambda , F^{101}) =
(\lambda^2 + \lambda - n \lambda - n^2) (\lambda^2 + \lambda)^{n - 1}$, and
 \\[0.7ex]
$(a3)$
$A(\lambda , F^{111}) = (\lambda + 1)^{2n - 1} (\lambda - 2n + 1)$.
\end{Theorem}

It turns out that if $F \in \cal{CF}$,
$z = +$, and $x, y \in \{0, +\}$, then
$A(\lambda, F^{xyz})$ is also uniquely defined by $v(F)$ and $c(F)$.
\begin{Theorem}
 \label{Fz+xy0+}
Let
$F  \in \cal{CF}$, $v(F) = n$, and $c(F) =  k$.
 Then
  \\[0.7ex]
$(a1)$
$A(\lambda , F^{00+}) = A(\lambda ^2, F) =  \lambda ^{2(n - k)} (\lambda ^{2k} - 1)$.
 \\[0.7ex]
$(a2)$
$A(\lambda , F^{0++}) = A(\lambda , F^{+0+}) =
(\lambda +1)^n~A(\frac{\lambda ^2}{\lambda +1},F)
=  \lambda ^{2(n-k)}[\lambda ^{2k} - (\lambda +1)^k]$, and
 \\[0.7ex]
$(a3)$
$A(\lambda , F^{+++}) = A(x_1, F) ~
 A(x_2, F)
= \lambda^{2n - 2k} (x_1^k - 1) (x_2^k - 1)$, where
\\[0.5ex]
$~~~~~~~x_1 = \frac{1}{2}(2\lambda +1 + \sqrt{4\lambda +1})$
  and $x_2 =
\frac{1}{2}(2\lambda +1 - \sqrt{4\lambda +1})$.
\end{Theorem}
\bp
By Remark \ref{subdividing},
$F^{00+}$ is the subdivision digraph of $D$. Obviously,
$F^{00+}$ is also a digraph-function with
$v(F^{00+}) = 2 v(F)$ and
$c(F^{00+}) =2 c(F)$. Therefore claim $(a1)$ follows from Theorem \ref{A(x,F)}.
\\[1ex]
\indent
${\bf (p1)}$
We prove claim $(a2)$.
 By Lemma \ref{Inverse}, it is sufficient to prove our claim when $F$ is  the inverse of a connected digraph-function.

Let $A = A(F), T= T(F)$ and $H = H(F)$.
We assume that $V(F)$ and $E(F)$ are  ordered in such a way that $T = A$ and  $H = I_n$. Since $F$ is the inverse of a degree-function, such orderings exist by Lemma \ref{crt} $(a4)$.
By Lemma \ref{AL} $(a2)$,
$A(F^l) = H^{\top}T$, and so $A(F^l) = A$.
Therefore

\begin{eqnarray}
\label{F0++}
 A(\lambda, F^{0++}) = \left|
                            \begin{array}{cc}
                              \lambda I_n &~~ -A                              \\[0.7ex]
                              -I_n &~~ \lambda I_n - A \\
                            \end{array}
                          \right|
 \end{eqnarray}

                                                  and
\begin{eqnarray}
\label{F+0+}
  A(\lambda, F^{+0+}) =
  \left|
                            \begin{array}{cc}
                              \lambda I_n - A &~~ -A                              \\[0.7ex]
                              -I_n &~~ \lambda I_n
                               \\
                            \end{array}
                          \right|.
   \end{eqnarray}

    By Lemma \ref{lemABCD}, we have from
 (\ref{F0++}) and (\ref{F+0+}):
\begin{eqnarray*}
A(\lambda, F^{0++}) = A(\lambda, F^{+0+})&=&
  |\lambda I_n|~|\lambda I_n - A -
  \lambda ^{-1} I_n A |
  \\
  &=&
  |\lambda ^2 I_n - (\lambda +1) A| =
  (\lambda +1)^n
  A(\frac{\lambda ^2}{\lambda +1},F).
\end{eqnarray*}

\noindent
Now the last equality follows from
Theorem \ref{A(x,F)} $(a1)$.

\vskip 1ex
${\bf (p2)}$
Finally we prove $(a3)$. By Lemma \ref{Inverse}, it is sufficient to prove our claim when $F$ is a digraph-function.
Let $A^l =  A(F^l)$.
By Lemma \ref{crt}, $F$ and $F^l$ are isomorphic, and so $n = v(F) = e(F) = v(F^l)$ and
$A(\lambda , F) = A(\lambda , F^l)$.
We assume that $V(F)$ and $E(F)$ are ordered in such a way that $T = I_n$ and, accordingly, $H^\top = A = A^l$. Since $F$ is a  digraph-function, by Lemma \ref{crt}  $(a4)$, such orderings of $V(F)$ and $E(F)$ exist.
By definition of $F^{+++}$,
\begin{eqnarray*}
  A(\lambda, F^{+++}) = \left|
                            \begin{array}{cc}
                              \lambda I_n - A&~~  - T
                              \\[0.7ex]
                             - H^\top &~~ \lambda I_n - A^l \\
                            \end{array}
                          \right| =
\left|
                            \begin{array}{cc}
                             \lambda I_n -A &~~ - I_n
                             \\[0.7ex]
                             - A  &~~ \lambda I_n - A \\
                            \end{array}
                          \right|.
\end{eqnarray*}
Adding the second row of the above block matrix to the first row, we obtain
\begin{eqnarray}
\label{D+++}
  A(\lambda, F^{+++}) = \left|
                            \begin{array}{cc}
                             \lambda I_n  &~~ -(\lambda + 1)I_n + A
                             \\[0.7ex]
                             -A  &~~ \lambda I_n - A \\
                            \end{array}
                          \right|.
\end{eqnarray}
By Lemma \ref{lemABCD}, we have from
(\ref{D+++}):
\begin{eqnarray*}
A(\lambda, F^{+++}) = |\lambda  I_n|~ |\lambda  I_n - A  - \lambda ^{-1} A ((\lambda + 1)I_n - A)| =
|\lambda ^2 I_n - (2\lambda + 1) A + A^2|.
\end{eqnarray*}
Therefore
$A(\lambda, F^{+++}) = A (x_1, F) ~A(x_2, F)$,
where $x_1$ and $x_2$ are the roots of
\\[0.7ex]
$x^2 -(2\lambda +1) x + \lambda ^2$, i.e.
$x_1 = \frac{1}{2}(2\lambda +1 + \sqrt{4\lambda +1})$
 and
$x_2 = \frac{1}{2}(2\lambda +1- \sqrt{4\lambda +1})$.
\ep

\subsection{Isomorphic and non-isomorphic  $xyz$-transformations of digraph-functions}
\label{isomorph-df-transform}

\indent

First we describe some pairs of triples
($xyz$, $x'y'z'$) with
$x,y,z, x', y', z' \in \{0,1,+,-\}$ such that
$D^{xyz}$ and $D^{x'y'z'}$ are isomorphic for every digraph-function and its inverse.
\begin{Theorem}
\label{F(xyz)-F(yxz)}
Let $F \in {\cal F}$.
Then
$F^{xyz}$ and $F^{yxz}$ are isomorphic {\em (and therefore, cospectral)}
 for all $x,y \in \{0, 1, +, -\}$ and $z \in \{0,1\}$.
\end{Theorem}

\begin{figure}[ht]
\label{Fig:r=1isomorph}
\begin{center}
\scalebox{0.2}[.2]{\includegraphics{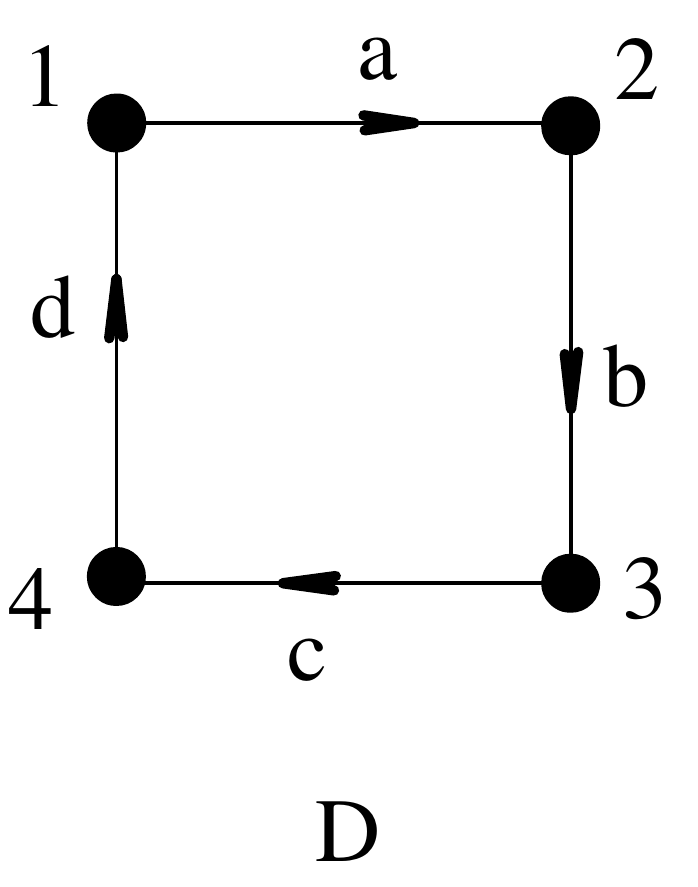}} \ \qquad
\scalebox{0.4}[.4]{\includegraphics{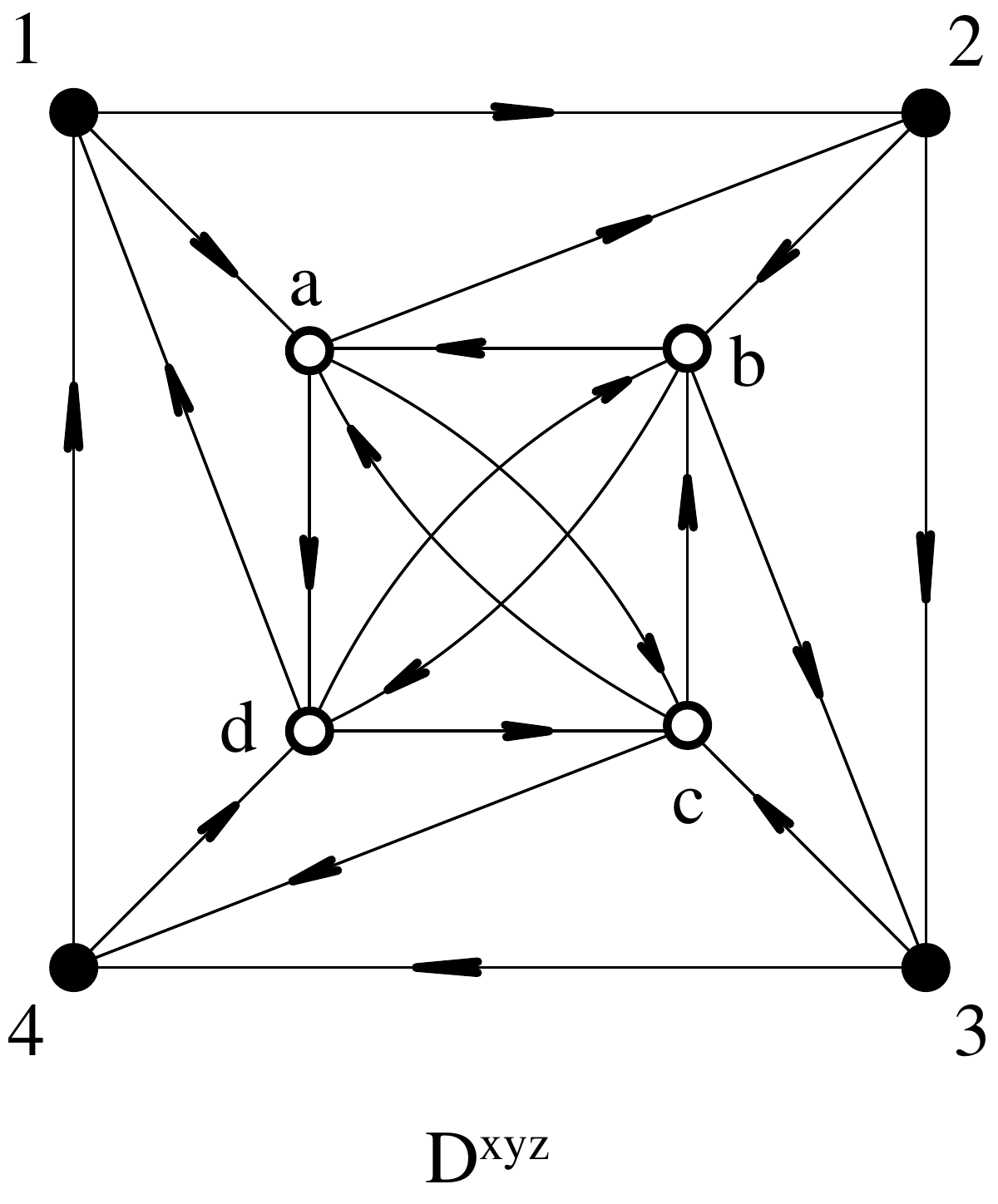}} \ \qquad
\scalebox{0.4}[.4]{\includegraphics{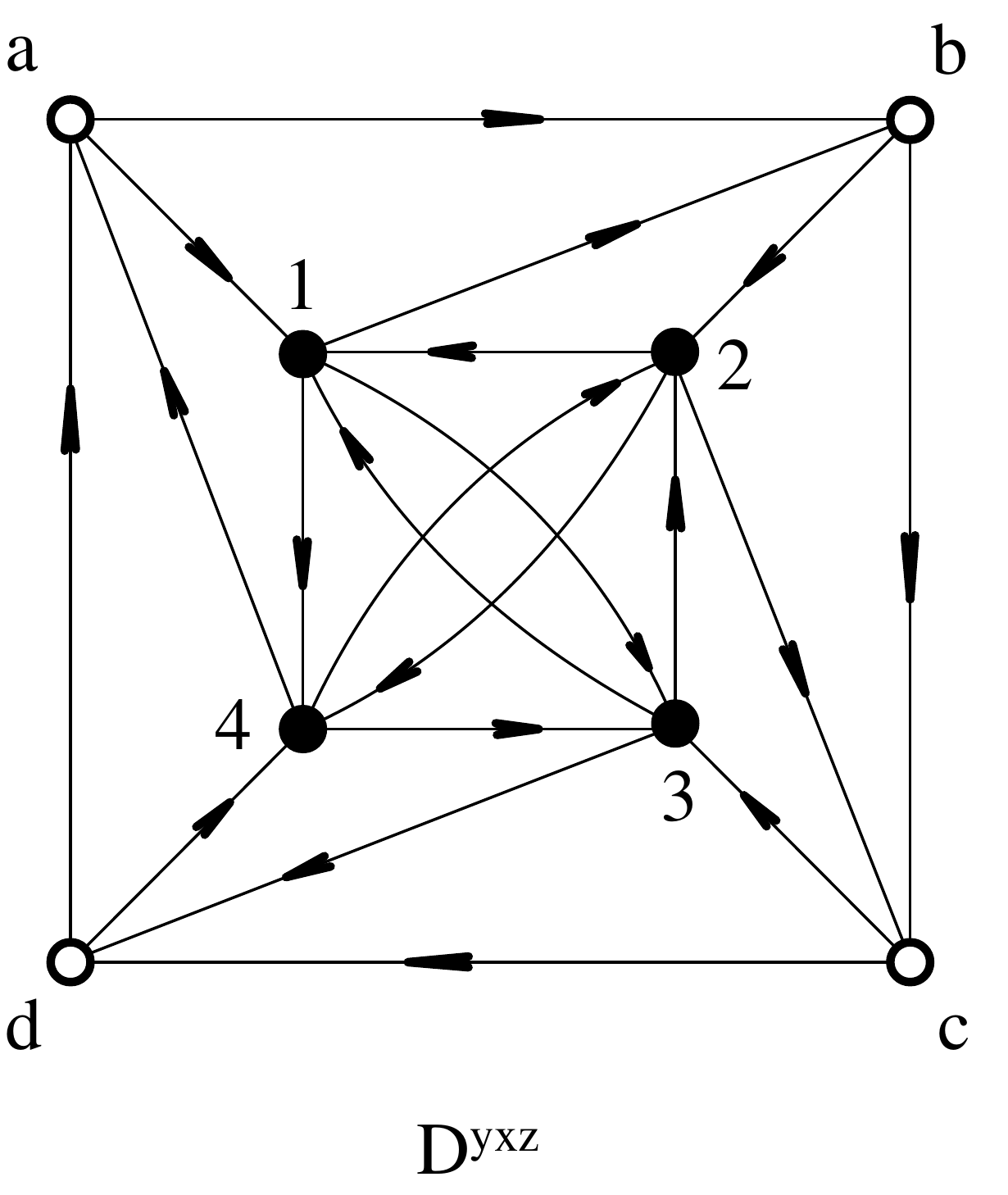}}
\end{center}
\caption{$D^{yxz} = \varepsilon [D^{xyz}]$, where $x =+, y= -, z= +$.}
\label{Dxyz=Dyxz}
\end{figure}

\bp We prove for $z =1$. Since $F \in {\cal F}$, by Lemma \ref{crt},
there is an isomorphism
$\alpha $ from $F$ to $F^l$.
We define a function $\varepsilon$ from
$V(F^{xy1})$ to $V(F^{yx1})$ such that $\varepsilon(w) = \alpha (w)$ if $w\in V(F)\subseteq V(F^{xy1})$ and
$\varepsilon(w) = \alpha ^{-1}(w)$ if $w\in E(F) \subseteq V(F^{xy1})$.
Then $\varepsilon $ is an isomorphism from
$F^{xy1}$ to $F^{yx1}$.
The proof for $z = 0$ is similar.
\ep
\\[1ex]
\indent
Our next result is on pairs of triples
($xyz$, $x'y'z'$) with
$x,y,z, x', y', z' \in \{0,1,+,-\}$ such that
$D^{xyz}$ and $D^{x'y'z'}$ are isomorphic for every 1-regular digraph.
\begin{Theorem}
\label{r=1isomorphism}
Let $D = (V, E)$ be a 1-regular digraph and  $x, y, z \in \{0,1,+,-\}$.
Given $w \in V \cup E$ let $\varepsilon (w) = t^{-1}(w)$ if
$w \in V$ and $\varepsilon (w) = h(w)$ if $w \in E$.
Then  $\varepsilon $ is an isomorphism from
$D^{xyz}$ to $D^{yxz}$.
\end{Theorem}
\bp
For $X \subseteq V\cup E$, we put
$\varepsilon [X] = \sum \{\varepsilon (x): x \in X\}$.
Clearly, $D^{xyz}$ is an edge disjoint union: $D^{xyz}
 = D^{xy0} \cup D^{00z}$.
Similarly, $D^{yxz}
= D^{yx0} \cup D^{00z}$.
If $x = y$, then clearly $D^{xy0} = D^{yx0}$.
Therefore we assume that $x \ne y$.
\\[1ex]
${\bf (p1)}$
Clearly, $\varepsilon $ is a bijection from $V\cup E$ to $V \cup E$, where
$V \cap E = \emptyset $.
Moreover, we have
\\[1ex]
{\sc Claim 1.}
\\[0.3ex]
{\em
$(a1)$
$\varepsilon |_V$ is an isomorphism from $D$ to $D^l$, and therefore also from $D^c$ to $(D^l)^c$ and
\\[0.5ex]
$(a2)$
$\varepsilon |_E$ is an isomorphism from $D^l$ to $D$, and therefore also from
$(D^l)^c$ to $(D)^c$.}
\\[1ex]
\indent
We also need the following fact.
\\[1ex]
{\sc Claim 2.}
{\em $\varepsilon $ is an automorphism of
$D^{00+}$}.
\\[1ex]
{\em Proof.} By Definition \ref{definition},
$D^{00+} = {\cal T}(D) \cup  {\cal H}(D)$.
By Remark \ref{subdividing}, $D^{00+}$ can be obtained from $D$ by subdividing each arc $e$ of $D$ into two arcs with a new vertex with label $e$. Therefore each component $\dot{C}$ of
 $D^{00+}$ is an even  directed cycle obtained from a directed cycle $C$ of $D$ by the above described subdivision. Then
 $\varepsilon$, restricted on $V(\dot{C})$, is an automorphism of $\dot{C}$, namely, a one-step rotation. Therefore  $\varepsilon $ is an automorphism of $D^{00+}$.
\epcl
\\[1ex]
${\bf (p2)}$
We prove our claim for $\{x,y\} = \{+,-\}$.
The case when $\{x,y\} \ne \{+,-\}$ can be proved similarly.
Obviously,
$D^{+-0} = D \cup (D^l)^c$ and
$D^{-+0} = D^c \cup (D^l)$. Therefore by {\sc Claim 1},
$\varepsilon $ is an isomorphism from
$D^{+-0}$ to $D^{-+0}$.
Hence our claim is true for $z = 0$.
By {\sc Claim 2}, $\varepsilon $ is also an automorphism of $D^{00+}$.
Therefore  $\varepsilon $ is an isomorphism from $D^{+-+} = D^{+-0} \cup D^{00+}$ to
$D^{-++} = D^{-+0} \cup D^{00+} $.
Thus, our claim is true for $z = +$.
Since $D^{-+1} = (D^{+-0})^c $ and
$D^{+-1} = (D^{-+0})^c$, our claim is also true for $z = 1$.
Also  $D^{+--} = (D^{-++})^c$ and
$D^{-+-} = (D^{+-+})^c$. Therefore our claim is also true for  $z = -$.
\ep
\\[1.5ex]
\indent
An example illustrating Theorem \ref{r=1isomorphism} is shown on Fig. \ref{Dxyz=Dyxz}.
%
\\[1ex]
\indent
Now we will describe some
pairs of triples
($xyz$, $x'y'z'$) with
$x,y,z, x', y', z' \in \{0,1,+,-\}$ such that
$D^{xyz}$ and $D^{x'y'z'}$ are not isomorphic for every non-regular connected digraph-function and its inverse.
\begin{Theorem}
\label{F(xyz)not-ismrF(yxz)}
Let $F \in \cal{CF}$ and $x,y, z \in \{0, 1, +, -\}$.
Then the following are equivalent:
\\[0.7ex]
$(a1)$
$F^{xyz}$ and $F^{yxz}$ are not isomorphic
and
\\[0.7ex]
$(a2)$
$z \in \{+,-\}$, $x \ne y$, and $F$ is not regular.
\end{Theorem}

\bp~
It is sufficient to prove our claim for connected digraph-functions.
By Theorems \ref{F(xyz)-F(yxz)} and
\ref{r=1isomorphism}, $(a1) \Rightarrow (a2)$.
We prove $(a2) \Rightarrow (a1)$.
By Lemma \ref{Gxyz} $(a1)$, it is sufficient to prove our claim for $z = +$.
Therefore we have to prove our claim for every two elements subsets $\{x,y\}$
of set $\{0, 1, +, -\}$, and so we have six corresponding cases to consider.

\begin{figure}[ht]
\begin{center}
\scalebox{0.7}[.7]{\includegraphics{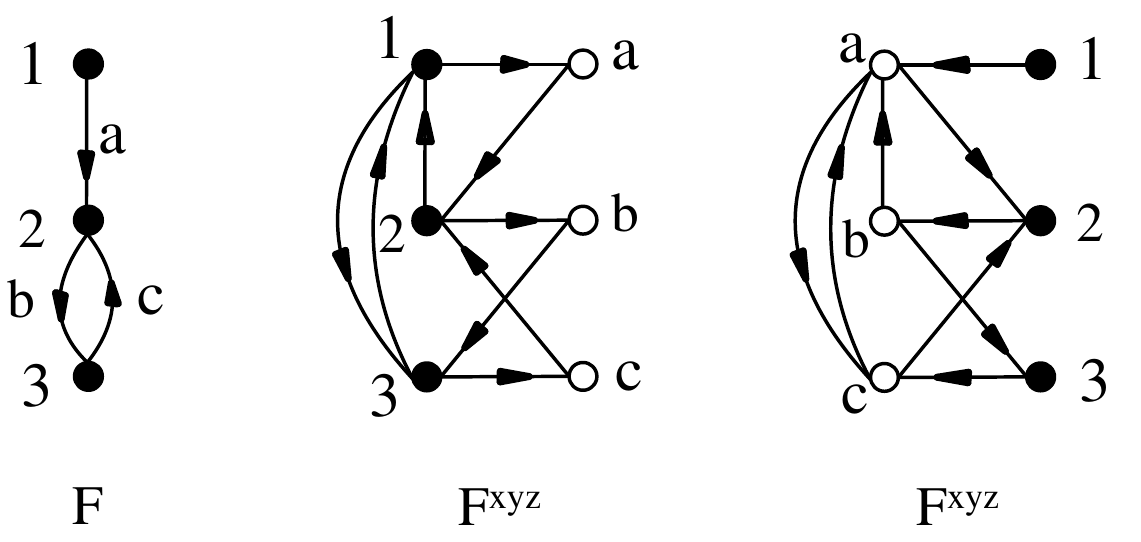}}
\end{center}
\caption
{$F^{xyz}$ and $F^{yxz}$ are not isomorphic, where $D\in \cal{CF}$, $x=-,~ y = 0$, and
$z = +$.}
\label{-0+*0-+}
\end{figure}

Let $V_s(D)$ denote the set of vertices of in-degree $s$ in a digraph $D$.
Obviously, $V_o(F^-) = V_o(F^1) =\emptyset$.
Since $F$ is not regular digraph-function, $0 < |V_o(F)| <  |V(F)|$. Since $z = +$, by Remark \ref{subdividing},
$V_o({\cal T}(F)\cup {\cal H}(F)) = V_o(F)$, and so also  $V_o(F \cup {\cal T}(F)\cup {\cal H}(F)) = V_o(F)$. Therefore by Definition
\ref{definition}, if $x \in \{+,0\}$, then
$V_o(F^{xy+}) = V_o(F) \ne \emptyset $.
From the above observations it follows that our claim is true for every two element subset  $\{x,y\}$ of set $\{0, 1, +, -\}$ distinct from $\{+, 0\}$ and $\{-,1\} $. Let $|V(F)| = r$. Then
  $|V_r(F^{1-+})| < |V_r(F^{-1+})| = r$ and so
 $F^{1-+}$ and $F^{-1+}$ are not isomorphic.
For $\{ x, y \} = \{ +, 0\}$, we have $d_{out}(v, F^{+0+}) = 2$ for each $v \in V_o(F^{+0+})$ but  $d_{out}(v, F^{0++}) = 1$ for each $v \in V_o(F^{0++})$. Therefore,   $F^{+0+}$ and $F^{0++}$ are not isomorphic.
\ep
\\[1.5ex]
\indent
An example illustrating Theorem  \ref{F(xyz)not-ismrF(yxz)} is shown in Fig. \ref{-0+*0-+}.

\subsection{Cospectral $xyz$-transformations of digraph-functions}
\label{cospectr-df-transform}

\indent

In subsection \ref{df-spectra} we described the adjacency polynomials of some
$xyz$-transformations of every degree-function and its inverse. As a byproduct, we have the following theorem on pairs of non-isomorphic digraphs $D$ and $F$ such that $D, F \in  \cal{CF}$ and
$D^{xyz}$ and $F^{xyz}$ are cospectral for some triples $xyz$ with $x,y,z \in \{0,1,+,-\}$.
\begin{Theorem}
\label{CospctrTransfD}
Let $D, ~F  \in \cal{F}$
and $x, y, z \in \{0, 1, +, -\}$.
If $D \sim^A F$, then
\\[1ex]
$(a1)$
$D^{xy0} \sim^A F^{xy0}$,
\\[1ex]
$(a2)$
$D^{00z} \sim^A F^{00z}$ for $z \neq -$,
 \\[1ex]
$(a3)$
$D^{0yz} \sim^A F^{0yz}$  for $y = z \in \{1, +\}$,
 \\[1ex]
$(a4)$
 $D^{x0z} \sim^A F^{x0z}$ for $x = z \in \{1, +\}$,
 \\[1ex]
$(a5)$
$D^{xyz} \sim^A F^{xyz}$  for $x = y = z \in \{1, +\}$,
and
\\[1ex]
$(a6)$ if $D$ and $F$ are not isomorphic, then the  above ${x'y'z'}$-transformations $D^{x'y'z'}$ and $F^{x'y'z'}$
of $D$ and $F$ in each of claims $(a1) - (a5)$ are also not isomorphic.
\end{Theorem}

\bp~ It is easy to prove $(a6)$.
Claims $(a1) - (a5)$ follow from Theorems \ref{Fz=0} - \ref{Fz+xy0+}.
\ep
\\[2ex]
\indent
Theorem \ref{F(xyz)not-ismrF(yxz)} provides the
characterization of all triples $xyz$ and $F \in \cal{CF}$ such that
$F^{xyz}$ and $F^{yxz}$ are not isomorphic.
In this subsection we provide the characterization of triples $xyz$ such that
$F^{xyz}$ and $F^{yxz}$ are cospectral for every $F \in \cal{F}$.
\begin{Theorem}
\label{F(xyz)--F(yxz)}
 Let $F$ be a
digraph-function or its inverse and $x,y,z \in \{0,1, +, -\}$.
Then
$F^{xyz}$  and $F^{yxz}$ are cospectral.
 \end{Theorem}

Obviously, if $F^{xyz}$  and $F^{yxz}$ are isomorphic, then they are  cospectral.
Therefore it is sufficient to prove the above theorem for $xyz$ and $F \in \cal{F}$ such that $F^{xyz}$  and $F^{yxz}$ are not isomorphic, i.e.  for $xyz$ and $F \in \cal{F}$ satisfying condition $(a2)$ in Theorem
\ref{F(xyz)not-ismrF(yxz)}: $F$ is not regular, $x \ne y$, and $z \in \{+, -\}$.

By Lemma \ref{Gxyz} $(a1)$, if
$D$ is a simple digraph and $x,y,z \in
\{0, 1, +, -\}$, then
digraphs $D^{xyz}$ and $D^{x'y'z'}$ are complement if and only if each of $\{x,x'\}$, $\{y,y'\}$, $\{z,z'\}$  is either
$\{0,1\}$ or  $\{+,-\}$.
Therefore by Corollary \ref{reciprocity}, it is sufficient to prove our theorem for $z = +$.
Thus, as in the proof of Theorem \ref{F(xyz)not-ismrF(yxz)}, there are  six two element subsets $\{x,y\}$ of set
$ \{0,1, +, -\}$  to consider.
Now Theorem  \ref{F(xyz)--F(yxz)} will follow from Lemmas \ref{F(+-+)--F(-++)},
\ref{x-in01,y-in+-}, and \ref{x=0,y=1}
below.
%
\begin{lemma}
 \label{F(+-+)--F(-++)}
 Let $F \in {\cal F}$. Then
$F^{+-+} \sim^AF^{-++}$.
\end{lemma}

\bp
 It is sufficient to prove our claim for a connected digraph-function.
Let $A = A(F)$, $A^l = A(F^l)$, $T= T(F)$, and $H = H(F)$.
We assume that $V(F)$ and $E(F)$ are  ordered in such a way that $T = I_n$ and  $H ^\top  = A$. Since $F$ is  a digraph-function, by Lemma \ref{crt} $(a4)$, such ordering exists.
By Lemma \ref{AL} $(a2)$,
$A^l = H^{\top}T = A$.

Therefore by definition of $F^{-++}$,
\begin{eqnarray}
\label{D-++}
  A(\lambda, F^{-++}) = \left|
                            \begin{array}{cc}
                             (\lambda +1) I_n - J_{nn} + A & - T
                              \\[0.7ex]
                              - H^\top &  \lambda I_n - A^l \\
                            \end{array}
                          \right| =
\left|
                            \begin{array}{cc}
                            (\lambda +1) I_n - J_{nn} + A  &  - I_n
                             \\[0.7ex]
                            - A   & \lambda I_n - A  \\
                            \end{array}
                          \right|.
\end{eqnarray}

Similarly,
by definition of $F^{+-+}$,
\begin{eqnarray}
\label{D+-+}
  A(\lambda, F^{+-+})
= \left|
                            \begin{array}{cc}
                               \lambda I_n - A & - T
                              \\[0.7ex]
                              - H^\top  & (\lambda +1) I_n - J_{nn} + A \\
                            \end{array}
                          \right|
=
\left|
                            \begin{array}{cc}
                            \lambda  I_n - A  & -I_n
                             \\[0.7ex]
                             - A   & (\lambda +1) I_n - J_{nn} + A \\
                            \end{array}
                          \right|.
\end{eqnarray}

Subtracting in the last matrix in (\ref{D-++})
the second block column from the first one, we obtain:
\begin{eqnarray}
\label{D-++S}
  A(\lambda, F^{-++})
=
\left|
                            \begin{array}{cc}
                            (\lambda +2) I_n - J_{nn} + A  &  - I_n
                             \\[0.7ex]
                            -  \lambda I_n   & \lambda I_n - A  \\
                            \end{array}
                          \right|.
\end{eqnarray}

Similarly subtracting in the last matrix in (\ref{D+-+})
the first block row from the second one, we obtain:
\begin{eqnarray}
\label{D+-+S}
  A(\lambda,F^{+-+})
=
\left|
                            \begin{array}{cc}
                            \lambda  I_n - A  & -I_n
                             \\[0.7ex]
                             -  \lambda  I_n   & (\lambda +2) I_n - J_{nn} + A \\
                            \end{array}
                          \right|.
\end{eqnarray}

Obviously, it is sufficient to prove our equality when
$\lambda $ is not an eigenvalue of $A$.
Then by Lemma \ref{lemABCD}, we have
from (\ref{D-++S}):
\begin{eqnarray}
\label{eq(-++)S}
 A(\lambda, F^{-++}) =
 |\lambda  I_n - A|~|(\lambda +2) I_n  - J_{nn} + A -\lambda I_n (\lambda  I_n - A)^{-1} I_n|,
  \end{eqnarray}
and \begin{eqnarray}
\label{eq(+-+)S}
 A(\lambda, F^{+-+}) =
 |\lambda  I_n - A|~
 |(\lambda +2) I_n -  J_{nn} + A -
 I_n (\lambda  I_n - A)^{-1} \lambda I_n|.
  \end{eqnarray}
Now we have  from
(\ref{eq(-++)S}) and (\ref{eq(+-+)S}):
$  A(\lambda, F^{-++}) =   A(\lambda, F^{+-+})$.
\ep
\begin{lemma}
\label{x-in01,y-in+-}
Let $F \in {\cal F}$.
Then $F^{xy+}\sim^A F^{yx+}$ for $x \in \{0,1\}$ and $y \in \{+,-\}$.
\end{lemma}

\bp The proof can be obtained from the proof of Lemma \ref{F(+-+)--F(-++)} as follows. Consider the last matrices $M_1$ and $M_2$ in
\eqref{D-++} and  \eqref{D+-+} in the proof of Lemma \ref{F(+-+)--F(-++)}.

Suppose that $y = +$.
Let us replace
$(\lambda +1) I_n - J_{nn} + A$ in $M_1$ and $M_2$ by $(\lambda + 1)I_n - J_{nn}$ if $x = 1$ and by $\lambda I_n$ if $x = 0$. Then we obtain the proofs of $F^{1-+} \sim^A F^{-1+}$ and $F^{0-+} \sim^A F^{-0+}$.

 Now suppose that $y = -$.
Let us replace $\lambda  I_n -  A$ in $M_1$ and $M_2$
by $(\lambda + 1)I_n - J_{nn}$ if $x = 1$ and by $\lambda I_n$ if $x = 0$.
Then we obtain the proofs of $F^{1-+} \sim^A F^{-1+}$ and $F^{0-+} \sim^A F^{-0+}$.
\ep
\begin{lemma}
\label{x=0,y=1}
Let $F \in {\cal F}$.
Then
$F^{01+}\sim^A F^{10+}$.
\end{lemma}

\bp
 It is sufficient to prove our claim for a connected digraph-function.
By definition of $F^{01+}$,
$A(\lambda ,F^{01+})$ is obtained from
the last matrix in \eqref{D+-+} by replacing $\lambda I_n - A$ by $\lambda I_n$ and $(\lambda +1) I_n - J_{nn} + A$ by $(\lambda +1) I_n - J_{nn}$.
Similarly, by definition of  $F^{10+}$,
$A(\lambda ,F^{10+})$ is obtained from
the last matrix in \eqref{D+-+} by replacing $\lambda I_n - A$ by
$(\lambda +1) I_n - J_{nn}$
and $(\lambda +1) I_n - J_{nn} + A$ by
$\lambda I_n$.
Now applying the first alternative of Lemma \ref{lemABCD} to  the matrix of $A(\lambda ,F^{01+})$ and the second alternative of Lemma \ref{lemABCD} to  the matrix of $A(\lambda ,F^{10+})$, we obtain:
\\[1ex]
\indent \qquad
$A(\lambda ,F^{01+}) = |\lambda I_n|~|)\lambda +1)I_n - J_{nn} - \lambda ^{-n}AI_n|$, and
\\[1ex]
\indent \qquad
$A(\lambda ,F^{10+}) = |\lambda I_n|~|)\lambda +1)I_n - J_{nn} - \lambda ^{-n}I_nA|$.
\\[1ex]
Therefore
$A(\lambda ,F^{01+}) = A(\lambda ,F^{10+})$.
\ep
\\[1ex]
\indent
From Theorems \ref{F(xyz)not-ismrF(yxz)} and \ref{F(xyz)--F(yxz)}  we have the following result.
\begin{corollary}
\label{???}
Let $F$ be a
digraph-function or its inverse and $x,y,z \in \{0,1, +, -\}$.
Then
 $F^{xyz}$ and $F^{yxz}$ are non-isomorphic and cospectral if and only if $F$ is non-regular, $x \ne y$ and
$z \in \{+,-\}$.
\end{corollary}

For the case when $\{ x, y \} = \{+, 0\}$ we have more general results which are given in Corollaries  \ref{Sp(+00=0+0)},
\ref{Sp(+01=0+1)} and \ref{Sp(+0-=0+-)} corresponding to $z = 0, 1, -$, respectively. It is also easy to prove the following result
for general digraphs when $z=+$.
\begin{Theorem} {\em \cite{zhang}}
Let $F$ be a general digraph. Then $F^{+0+} \sim^A F^{0++}$.
\end{Theorem}

\section{More on cospectral  transformation
digraphs}
\label{cospectral}

\indent

The above results give various constructions providing adjacency cospectral digraphs.
Here is an overview of some of those constructions providing infinitely many pairs of cospectral and, obviously,
non-isomorphic digraphs.
\begin{Theorem}
\label{cospectral-constructions} Let $D$ and $F$ be digraphs. Then
\\[1ex]
$(a1)$ $D\sim ^A F \Rightarrow D^{-1}\sim ^A F^{-1}$ and $D^l\sim ^AF^l$,
 \\[1ex]
$(a2)$
 if $D$ and $F$ are simple regular digraphs, then
$D\sim ^A F \Rightarrow D^c \sim ^A F^c ~and~
D^{xyz}\sim ^A F^{xyz}$ for $x,y,z \in \{0,1,+,-\}$,
 \\[1ex]
$(a3)$
 if $D$ and $F$ are simple digraphs, then
$D\sim ^A F \Rightarrow D^{xy0}\sim ^A F^{xy0}$ for $x,y \in \{0,1,+,-\}$,
$D\sim ^A F \Rightarrow D^{+10}\sim ^A F^{+10},~and~
D\sim ^A F \Rightarrow D^{1+0}\sim ^A F^{1+0} $,
\\[1ex]
$(a4)$
 if $D$ and $F$ are digraphs, then
$D \sim ^A F \Rightarrow D^{+00}\sim ^A F^{+00},~D^{0+0}\sim ^A F^{0+0} ,~and~
D^{++0}\sim ^A F^{++0} $,
\\[1ex]
$(a5)$
 if $D$ and $F$ are general regular digraphs, then  $D^{xyz}\sim ^A F^{xyz}$ for $z = 1$,
$xy \in \{+-,0+,++\}$, and $D^{+01}\sim ^A D^{0+1}$,  $D^{+0-}\sim ^A D^{0+-}$,
\\[1ex]
$(a6)$ if $D, F \in \mathcal{CF}$, $v(D) = v(F)$ and
$c(D) = c(F)$, then $D \sim ^A F$,
and
\\[1ex]
$(a7)$
 if $D$ is a general  digraph,  then
 $D^{+00}\sim ^A D^{0+0}$.
 \end{Theorem}

Now  we  describe some more constructions that provide cospectral non-isomorphic and non-regular digraphs.

Let $D$ and $D'$ be disjoint digraphs,
$X \subseteq D$, $X' \subseteq D'$, $X \ne \emptyset$,
and $\pi $  a bijection from $X$ to $X'$.
Let $DX\pi X'D'$ denote the digraph obtained from $D$ and $D'$ by identifying vertex $x$ in $D$ with the vertex $\pi (x)$ in $D'$ for every $x \in X$.

Given a digraph $D$, let
$V_{in}(D) = \{v \in V(D): d_{out}(v) = 0\}$ and
$V_{out}(D) = \{v \in V(D): d_{in}(v) = 0\}$.
A digraph $D$ is called {\em acyclic} if $D$ has no directed cycles.
\\[1ex]
\indent
It is easy to prove the following:
\begin{Theorem}
\label{DpiD'}
Let $D$ and $D'$ be disjoint digraphs,
$F = DX\pi X'D'$, and $n = v(F)$, $k = v(D)$.
Suppose that $D'$ is an acyclic digraph and
$X' \subseteq V_{in}(D')$ or  $X' \subseteq V_{out}(D')$.
Then
$A(\lambda , F) = \lambda ^{n - k} A(\lambda , D)$.
\end{Theorem}

Obviously, Lemma \ref{A(x,F)} is a particular case of Theorem \ref{DpiD'}.

From Lemma \ref{PL} and Theorem \ref{DpiD'} we have:
\begin{Theorem}
\label{CospctrDl,D'}
Let
 $D$ and $D'$ be disjoint digraphs, $D'$ an acyclic digraph,
$F = DX\pi X'D'$, where
$X' \subseteq V_{in}(D)$ or
$X' \subseteq V_{out}(D)$, and $\pi $ is a bijection from $X$ to $X'$.
Suppose that $|V(D') \setminus X'| = e(D)  - v(D)$.
Then $A(\lambda , D^l) = A(\lambda , F)$.
\end{Theorem}

Theorems \ref{DpiD'} and \ref{CospctrDl,D'} give  constructions that provide an
infinite variety of
non-isomorphic cospectral digraphs.

\section{Some remarks}
\label{remarks}

\indent

${\bf (R1)}$
Notice that all the factors of the adjacency polynomials we present for $D^{xyz}$ ($x,y,z\in \{0,1,+,-\}$) are polynomials in $\lambda$ of degree one or two. Hence the explicit formula for the spectrum of $D^{xyz}$ can be given in terms of the spectrum of $D$, as in Corollaries \ref{Sp(+01=0+1)} and \ref{Sp(-01)}.
\\[1.5ex]
\indent
${\bf (R2)}$
Let ${\cal R}$ denote the set of simple regular digraphs.
Obviously, if $D \in {\cal R}$, then
$D^c \in {\cal R}$,
$D^{-1}  \in {\cal R}$, and $D^l \in {\cal R}$.
If $D$ is an $r$-regular digraph, then
$D^{+++}$ is $2r$-regular and $G^{---}$ is
$(v(D) + e(D) - 2r -1)$-regular,
and so if $D \in {\cal R}$, then $D^{+++}, D^{---} \in {\cal R}$.
In other words, the set ${\cal R}$ of simple regular digraphs is closed under  the $(-1)$-operation of taking the inverse,
{\em $c$-operation}, {\em $l$-operation},
{\em $(+++)$-operation},
and {\em $(---)$-operation}.
Therefore using the corresponding results described above, one can give an algorithm (and the computer program) that for
any series $Z$ of $(-1)$-, $c$-, $l$-, $(+++)$-, and $(---)$-operations and the spectrum $S_a(D)$ of any
$r$-regular digraph $D$ provides the formula of the spectrum of digraph $F$ obtained from $D$ by the  series $Z$ of operations in terms of $r$,
$v(D)$, and $S_a(D)$.
\\[1.5ex]
\indent
${\bf (R3)}$
Suppose that a regular digraph $D$ is uniquely defined by its adjacency  spectrum. Does it necessarily mean that $D^{xyz}$ is also uniquely defined by its  adjacency   spectrum for all or for some
$x,y,z \in \{ +, -\}$ ?
\\[1.5ex]
\indent
${\bf (R4)}$
Obviously, if regular digraphs $D$ and $F$ are isomorphic, then $D^{xyz}$ and $F^{xyz}$ are also isomorphic.
A natural question is whether
there exist two non-isomorphic regular
digraphs $D$ and $F$  such that  $D^{xyz}$ and $F^{xyz}$ are isomorphic for some
$x,y,z \in \{+, -\}$ ?
\\[1.5ex]
\indent
${\bf (R5)}$
Here is another definition of the digraph ${xyz}$-transformations for $x,y,z \in \{0,1,+,-\}$ that is valid for digraphs with loops but without multiple arcs. This definition is using the notion of  $K_\circ $-complement of $D$.

Recall that $K_{\circ} = (V, E)$, where
$E = V\times V$ is a {\em complete digraph}, and so every vertex in $K_{\circ}$ has a loop.
Given a  digraph $D$ with $V = V(D) = V(K_\circ )$ and $E(D) \subseteq V(K_\circ )$, let
$D^c_\circ = K_\circ \setminus E(D)$. Digraph
 $D^c_\circ $ is called the {\em $K_\circ $-complement of} $D$.

\begin{definition}
\label{definition2}
Given a digraph $D$ and three variables
$x, y, z \in \{0,1, +, -\}$, the {\em ${xyz}$-transformation $D^{xyz}_\circ$  of} $D$ is the digraph
such that
$D^{xy0}_\circ= D^x_\circ \cup [D^l]^y_\circ$
and $D^{xyz}_\circ = D^{xy0}_\circ \cup W$, where
$W = {\cal T}(D) \cup {\cal H}(D)$  if $z = +$,
$W = {\cal T}^c(D) \cup {\cal H}^c(D)$ if $z = -$, and
$W$ is the union of complete $(V,E)$-bipartite and $(E,V)$-bipartite digraphs if $z = 1$.
\end{definition}

Obviously, if $x,y \in \{0,+\}$, then all above results
for general digraphs are also valid  for ${xyz}$-transformation $D^{xyz}_\circ$  of a digraph $D$.
Here is the analog of Theorem \ref{-01} for
$D^{xyz}_\circ$,
i.e. when
$A(D^c_\circ) = J_{nn} - A(D)$ and
$A([D^l]^c_\circ) = J_{mm} - A(D^l)$.
\begin{Theorem}
\label{-01K0}
Let $D$ be an $r$-regular digraph with $n$ vertices and $m$ arcs.
  Then
 \[ A(\lambda , D^{-01}_\circ) = (-1)^n \lambda ^{m-1}
  (\lambda ^2  + r\lambda  - n\lambda  - mn) (\lambda  + r)^{-1}A(- \lambda , D).\]
   \end{Theorem}

 The proof of this theorem is similar to the proof of Theorem \ref{xy0+r}.
\begin{corollary}
\label{Sp(-01)K0}
Let $D$ be an $r$-regular digraph with $n$ vertices and $m$ edges and
\\
$S_a(D) = \{\alpha _i: i = 1, \ldots, n\}$, where $\alpha _n = r$.
Then $D^{-01}_\circ $ has $m-1$ eigenvalues equal to zero,  $n-1$  eigenvalues  $- \alpha_i $, where $i = 1, \ldots , n-1$, and two additional eigenvalues
\[ \frac{1}{2}( n - r  \pm \sqrt{(n - r )^2 + 4 r n^2 }). \]
\end{corollary}



\begin{center}
{\LARGE Appendix }
\label{appendix2}
\end{center}

Let $D$ be an $r$-regular digraph with $n$ vertices and $m$ edges,
and so $m = nr$.
The tables below provide the formulas for
$A(\lambda, D^{xyz})$ for all $x, y, z \in \{0, 1, + , -\}$ in terms of $n$, $r$, $m$, and the adjacency polynomials of $D$.
\\[3ex]
\indent
{\sc The list of }
$A(\lambda, D^{xyz})$
{\sc with}
$z=0$.
\\[1.5ex]
\begin{tabular}
{|l|l|l|}
  \hline
  & $xyz$ &  $A(\lambda, D^{xyz})$  \\
 \hline
1 & $0~0~0$ & $\lambda^{m+n}$  \\
 \hline
2 & $1~0~0$ & $\lambda^{m}(\lambda - n + 1) (\lambda + 1)^{n-1}$   \\
 \hline
3 & $+0~0$ & $\lambda^m A(\lambda, D)$   \\
 \hline
4 & $-0~0$ & $(-1)^n
\lambda^{m}
(\lambda - n + r + 1)(\lambda + r + 1)^{-1}
A(-\lambda - 1, D)$ \\
 \hline
5 & $0~1~0$ & $\lambda^{n} (\lambda + 1)^{m - 1} (\lambda - m + 1) $   \\
 \hline
6 & $1~1~0$ & $ (\lambda + 1)^{m + n -2} (\lambda -n + 1)(\lambda - m + 1) $   \\
 \hline
7 & $+1~0$ & $ (\lambda + 1)^{m-1} (\lambda - m + 1) ~A(\lambda, D)$   \\
 \hline
8 & $-1~0$  & $(-1)^n (\lambda + 1)^{m-1} (\lambda -m + 1)
(\lambda - n + r + 1)(\lambda + r + 1)^{-1}
A(-\lambda - 1, D)$ \\
 \hline
9 & $0+0$  & $\lambda ^m  A(\lambda , D)$   \\
 \hline
10 & $1+0$  & $\lambda^{m - n} (\lambda + 1)^{n - 1} (\lambda - n + 1) A(\lambda , D)$  \\
 \hline
11 & $++0$  & $\lambda^{m - n} A(\lambda , D)^2$ \\
 \hline
12 & $-+0$  & $(-1)^n \lambda^{m - n}
(\lambda -n + r + 1)(\lambda + r + 1)^{-1}
A(- \lambda  - 1, D)~ A(\lambda , D)$  \\
 \hline
13 & $0-0$  & $(-\lambda)^n (\lambda + 1)^{m - n}
(\lambda - m + r + 1)(\lambda + r + 1)^{-1}
A(-\lambda - 1, D)$  \\
 \hline
14 & $1-0$  & $ (-1)^n (\lambda + 1)^{m -1}  (\lambda -n + 1)
(\lambda - m + r + 1)(\lambda + r + 1)^{-1}
 A(-\lambda - 1, D)$  \\
 \hline
15 & $+-0$  & $(-1)^n  (\lambda + 1)^{m-n}
(\lambda - m + r + 1)(\lambda + r + 1)^{-1}
 A( - \lambda - 1, D)
~A(\lambda , D)$   \\
 \hline
 16 & $--0$ &
  $( \lambda + 1)^{m - n}
 (\lambda -n + r + 1)(\lambda - m + r + 1)(\lambda + r + 1)^{-2}
%
%
A ( - \lambda - 1, D)^2 $   \\
 \hline
\end{tabular}
\\[3ex]

\newpage

{\sc The list of }
$A(\lambda, D^{xyz})$
{\sc with}
$z=1$.
\\[1.5ex]
\begin{tabular}
{|l|l|l|}
  \hline
  & $xyz$ &  $A(\lambda, D^{xyz})$  \\
 \hline
1 & $0~0~1$ & $\lambda^{m+n - 2} (\lambda^2 - m n)$  \\
 \hline
2 & $1~0~1$ & $\lambda^{m - 1} (\lambda + 1)^{n - 1} (\lambda^2 + \lambda - n \lambda - m n)$   \\
 \hline
3 & $+0~1$ & $\lambda^{ m - 1}
(\lambda^2 - r \lambda - m n)(\lambda - r)^{-1}
A(\lambda, D)$   \\
 \hline
4 & $-0~1$ & $(-1)^n
\lambda^{m - 1}
(\lambda^2 + \lambda + r\lambda - n \lambda - m n)(\lambda + r + 1)^{-1}
  A(-\lambda - 1, D)$ \\
 \hline
5 & $0~1~1$ & $\lambda^{n - 1} (\lambda + 1)^{m - 1} (\lambda ^2 + \lambda - m \lambda - mn) $   \\
 \hline
6 & $1~1~1$ & $ (\lambda + 1)^{m + n -1} (\lambda -m - n + 1) $   \\
 \hline
7 & $+1~1$ & $ (\lambda + 1)^{m-1}
((\lambda - r)(\lambda + 1 - m)- m n)(\lambda - r)^{-1}
A(\lambda, D)$   \\
 \hline
8 & $-1~1$  & $(-1)^n (\lambda + 1)^{m-1}
((\lambda + 1)(\lambda + r + 1 - m - n) - m r)(\lambda + r + 1)^{-1}
A(-\lambda - 1, D)$ \\
 \hline
9 & $0+1$  & $\lambda^{m - 1}
(\lambda^2 - r \lambda - m n)(\lambda - r)^{-1}
A(\lambda , D)$   \\
 \hline
10 & $1+1$  & $\lambda^{m - n} (\lambda + 1)^{n - 1}
((\lambda - r)(\lambda + 1 - n) - m n)(\lambda - r)^{-1}
A(\lambda , D)$  \\
 \hline
11 & $++1$  & $\lambda^{m - n}
((\lambda - r)^2 - m n)(\lambda - r)^{-2}
A(\lambda , D)^2$ \\
 \hline
12 & $-+1$  & $(-1)^n \lambda^{m - n}
((\lambda - r)(\lambda + r + 1 - n) - m n)((\lambda - r)(\lambda + r + 1))^{-1}A(- \lambda  - 1, D) $
\\
& &$A(\lambda , D)$  \\
 \hline
13 & $0-1$  & $(- 1)^n \lambda^{n-1} (\lambda + 1)^{m - n}
(\lambda (\lambda + r + 1 - m) - m n) (\lambda + r + 1)^{-1}
A(-\lambda - 1, D)$  \\
 \hline
14 & $1-1$  & $ (-1)^n (\lambda + 1)^{m -1}
((\lambda + 1)(\lambda + r + 1 - m - n) - n r)(\lambda + r + 1)^{-1}
 A(-\lambda - 1, D)$  \\
 \hline
15 & $+-1$  & $(-1)^n  (\lambda + 1)^{m-n}
 ((\lambda - r)(\lambda + r + 1 - m) - m n )((\lambda - r)(\lambda + r + 1))^{-1}$
\\
& & $ A( - \lambda - 1, D)
~A(\lambda , D)$   \\
 \hline
 16 & $--1$ & $( \lambda + 1)^{m - n}
(\lambda + r + 1 - m - n)(\lambda + r + 1)^{-1}
A ( - \lambda - 1, D)^2 $   \\
 \hline
\end{tabular}

\newpage

{\sc The list of }
$A(\lambda, D^{xyz})$
{\sc with}
$z=+$.
\\[1.5ex]
\begin{tabular}{|l|l|l|}
  \hline
    & $xyz$ &  $A(\lambda, D^{xyz})$  \\
 \hline
1 & $0~0~+$ & $\lambda^{m-n} A(\lambda^2, D)$  \\
 \hline
2 & $1~0~+$ & $\lambda^{m - n}
(\lambda^2 +\lambda - n \lambda - r)(\lambda^2 + \lambda - r)^{-1}
 A(\lambda^2 + \lambda, D)$   \\
 \hline
3 & $+0~+$ & $\lambda^{ m - n} (\lambda + 1)^n
A(\lambda^2(\lambda + 1)^{-1}, D)$   \\
 \hline
4 & $-0~+$ & $(1 - \lambda)^n
\lambda^{m - n}
(\lambda^2 + \lambda + r(\lambda -1) - n \lambda ) (\lambda^2 + \lambda + r(\lambda - 1))^{-1}
A((\lambda^2 + \lambda)(1- \lambda)^{-1}
, D)$ \\
 \hline
5 & $0~1~+$ & $ (\lambda + 1)^{m - n}
(\lambda^2 + \lambda - r - m \lambda)(\lambda^2 + \lambda - r)^{-1}
A(\lambda^2 + \lambda, D)$   \\
 \hline
6 & $1~1~+$ & $ (\lambda + 1)^{m - n}
((\lambda + 1)^2 -n(r+1)(\lambda +1) - r + m n)((\lambda + 1)^2 - r)^{-1}
A((\lambda+1)^2, D)$   \\
 \hline
7 & $+1~+$ & $ (\lambda + 1)^{m-n} (\lambda + 2)^n
(\lambda^2 + \lambda - r(\lambda + 2) - m(\lambda - r))(\lambda^2 -(r-1) \lambda
- 2r)^{-1}$
 \\
& & $ A((\lambda^2 + \lambda)(\lambda + 2)^{-1}
, D)$   \\
 \hline
8 & $-1~+$  & $(- \lambda)^n (\lambda + 1)^{m-n}
((\lambda + 1)^2 + r \lambda -n (\lambda + 1)- m(\lambda - n + r + 1))((\lambda + 1)^2 + r \lambda)^{-1}$
\\
& &$A(-(\lambda + 1)^2\lambda ^{-1}
, D)$ \\
 \hline
9 & $0++$  & $\lambda^{m - n} (\lambda +1)^n
A(\lambda^2(\lambda + 1)^{-1}
, D)$   \\
 \hline
10 & $1++$  & $\lambda^{m - n} (\lambda + 2)^{n}
(\lambda^2 + \lambda - r(\lambda + 2) + n (r - \lambda))(\lambda^2 + \lambda - r(\lambda + 2))^{-1}$
\\
& & $A((\lambda^2 + \lambda)(\lambda + 2)^{-1}
, D) $  \\
 \hline
11 & $+++$  & $\lambda^{m - n}  A(\frac{1}{2}(2\lambda + 1 + \sqrt{4\lambda + 1}) , D)  ~A(\frac{1}{2}(2\lambda + 1 - \sqrt{4\lambda + 1}) , D) $ \\
 \hline
12 & $-++$  & $(-1)^n \lambda^{m - n}
(\lambda^2 - (n-1) \lambda + m- 2r - r^2
) (\lambda^2 + \lambda - r^2 - 2r)^{-1}$
%
%
\\
& &$A(- 1 + \sqrt{\lambda^2 + \lambda + 1}, D) ~A(- 1 - \sqrt{\lambda^2 + \lambda + 1} , D)$  \\
 \hline
13 & $0-+$  & $( 1 -\lambda)^n (\lambda + 1)^{m - n}
(\lambda^2 +\lambda + r(\lambda - 1) - m \lambda)(\lambda^2 +\lambda + r(\lambda - 1))^{-1}$
\\
& & $A((\lambda^2 + \lambda)( 1 - \lambda)^{-1}
, D)$  \\
 \hline
14 & $1-+$  & $ (-\lambda)^n (\lambda + 1)^{m -n}
((\lambda + 1)^2 + r \lambda - n(\lambda + 1)- m (\lambda + 2 - n ))((\lambda + 1)^2 + r \lambda )^{-1}$
\\
& & $ A(-(\lambda + 1)^2 \lambda^{-1}
, D)$  \\
 \hline
15 & $+-+$  & $(-1)^n  (\lambda + 1)^{m-n}
(\lambda^2 - (n-1) \lambda + m - r^2 - 2r
) (\lambda^2 + \lambda - r^2 - 2r)^{-1}$
%
%
\\
& &$ A( -1 + \sqrt{\lambda^2 + \lambda + 1}, D)
~A(-1 - \sqrt{\lambda^2 + \lambda + 1}, D) $   \\
 \hline
 16 & $--+$ & $( \lambda + 1)^{m - n}
  ((\lambda + 1)^2 + r(2\lambda + 1) + r^2 - n (\lambda + 1) - m(\lambda + 2 - n + r))((\lambda + 1)^2 + $
 \\
& & $  r(2\lambda + 1) + r^2 )^{-1}
 A ( \frac{1}{2}(-2\lambda - 1 + \sqrt{-4\lambda - 3}), D)~
A ( \frac{1}{2}(-2\lambda - 1 - \sqrt{-4\lambda - 3}), D)$ \\
 \hline
\end{tabular}

\newpage

{\sc The list of }
$A(\lambda, D^{xyz})$
{\sc with}
$z = -$.
\\[1.5ex]
\begin{tabular}{|l|l|l|}
  \hline
  & $xyz$ &  $A(\lambda, D^{xyz})$  \\
 \hline
1 & $0~0~-$ & $\lambda^{m - n} (\lambda^2 - r + 2m - mn) (\lambda^2 - r)^{-1} A(\lambda^2, D)$  \\
 \hline
2 & $1~0~-$ & $\lambda^{m - n}
(\lambda^2 - (n-1) \lambda - r + 2m - mn)
(\lambda^2 + \lambda - r)^{-1}
A(\lambda^2 + \lambda, D)$   \\
 \hline
3 & $+0~-$ & $\lambda^{ m - n} (\lambda + 1)^n
(\lambda^2 - r\lambda - r + 2m - m n)(\lambda^2 - r\lambda - r)^{-1}
%
A(\lambda^2(\lambda + 1)^{-1}
, D)$   \\
 \hline
4 & $-0~-$ & $(1 - \lambda)^n \lambda^{m - n}
(\lambda^2 + \lambda (r + 1- n) - r  + 2m -
mn )
(\lambda^2 + ( r + 1)\lambda - r)^{-1}$
\\
& & $ A((\lambda^2 + \lambda)(1- \lambda)^{-1}
 , D)$ \\
 \hline
5 & $0~1~-$ & $ (\lambda + 1)^{m - n}
(\lambda^2 + \lambda - r - m (\lambda + n - 2))(\lambda^2 + \lambda - r)^{-1}
A(\lambda^2 + \lambda, D)$   \\
 \hline
6 & $1~1~-$ & $ (\lambda + 1)^{m - n}
((\lambda + 1)^2 - r -n(\lambda +1 + r \lambda - r))((\lambda + 1)^2 - r)^{-1}
 A((\lambda+1)^2, D)$   \\
 \hline
7 & $+1~-$ & $ (\lambda + 1)^{m-n} (\lambda + 2)^n
 (\lambda^2 + \lambda - r(\lambda + 2) - m(\lambda + n - r - 2))
 (\lambda^2 -(r-1) \lambda - 2r)^{-1}$
%
\\
& & $ A((\lambda^2 + \lambda)(\lambda + 2)^{-1}
, D)$   \\
 \hline
8 & $-1~-$  & $(- \lambda)^n (\lambda + 1)^{m-n}
((\lambda + 1)^2 + r \lambda -n (\lambda + 1 - r)- m(\lambda + r))((\lambda + 1)^2 + r \lambda)^{-1}$
\\
& & $A(-(\lambda + 1)^2\lambda ^{-1}
, D)$ \\
 \hline
9 & $0+-$  & $\lambda^{ m - n} (\lambda + 1)^n
(\lambda^2 - r\lambda + 2m - m n - r)(\lambda^2 - r\lambda -r )^{-1}
%
%
A(\lambda^2(\lambda + 1)^{-1}
, D)$   \\
 \hline
10 & $1+-$  &
$\lambda^{m - n} (\lambda + 2)^{n}
(\lambda^2  - (n + r -1)\lambda  -
mn+ 3m - 2)
(\lambda^2 - (r -1) \lambda - 2r)^{-1}
%
$
\\
& & $A((\lambda^2 + \lambda)(\lambda + 2)^{-1}
, D) $  \\
 \hline
11 & $++-$  &
$\lambda^{m - n}
(\lambda^2 - 2r \lambda - mn + 2m + r^2
- r)(\lambda^2 - 2 \lambda + r^2 - r)^{-1}$
%
%
\\
& & $A(\frac{1}{2}(2\lambda + 1 + \sqrt{4\lambda + 1}) , D) ~ A(\frac{1}{2}(2\lambda + 1 - \sqrt{4\lambda + 1}) , D) $ \\
 \hline
12 & $-+-$  & $(-1)^n \lambda^{m - n}
(\lambda^2  - (n -1) \lambda
 - m n + 3m - r^2 - 2r)(\lambda^2 + \lambda - r^2 - 2r)^{-1}$
%
%
\\
& & $A(- 1 + \sqrt{\lambda^2 + \lambda + 1}, D)~ A(- 1 - \sqrt{\lambda^2 + \lambda + 1} , D)$  \\
 \hline
13 & $0--$  & $( 1 -\lambda)^n (\lambda + 1)^{m - n}
(\lambda^2 + (r - m +1)\lambda + (r - m)
( n - 1)) (\lambda^2  + (r +1)\lambda - r)^{-1}$
%
\\
& & $A((\lambda^2 + \lambda)
(1 - \lambda)^{-1}
, D)$  \\
 \hline
14 & $1--$  & $ (-\lambda)^n (\lambda + 1)^{m -n}
((\lambda + 1)^2 + r \lambda -n(\lambda + 1+ r \lambda))((\lambda + 1)^2 + r \lambda )^{-1}$
\\
& & $ A( -(\lambda + 1)^2 \lambda^{-1}
 , D)$  \\
 \hline
15 & $+--$  & $(-1)^n  (\lambda + 1)^{m-n}
(\lambda^2 + \lambda -2r -r^2 + m(r + 2 - n - \lambda) )(\lambda^2 + \lambda -2r -r^2)^{-1}$
\\
& & $A( -1 + \sqrt{\lambda^2 + \lambda + 1}, D)
~A(-1 - \sqrt{\lambda^2 + \lambda + 1}, D) $   \\
 \hline
 16 & $---$ & $( \lambda + 1)^{m - n}
((\lambda + 1)^2 + r(2\lambda + 1) + r^2 - n (\lambda + 1 + r \lambda + r^2))((\lambda + 1)^2 +$
\\
 & &  $ r(2\lambda + 1) + r^2 )^{-1}
A ( \frac{1}{2}(-2\lambda - 1 + \sqrt{-4\lambda - 3}), D)
~A ( \frac{1}{2}(-2\lambda - 1 - \sqrt{-4\lambda - 3}), D) $   \\
 \hline
\end{tabular}
\end{document}